\documentclass[a4,11pt]{article} 
\usepackage{mathrsfs}
\usepackage{adjustbox}
\usepackage{amsthm} 
\usepackage{mathrsfs}
\usepackage{url} 
\usepackage{hyperref} 
\usepackage{amsmath, amssymb}
\usepackage{graphicx}
\usepackage{color}
\usepackage{bm}
\usepackage[textwidth=14.17cm]{geometry}

\newtheorem{theorem}{Theorem}
\newtheorem{proposition}[theorem]{Proposition}%
\newtheorem{lemma}[theorem]{Lemma}
\newtheorem{corollary}[theorem]{Corollary}

\newtheorem{definition}{Definition}%

\usepackage{amsmath, amssymb, graphicx}
\usepackage{mathrsfs}
\usepackage{url} 
\usepackage{hyperref} 
\usepackage{color}
\usepackage{bm}
\usepackage[margin=0pt]{subfig} 

\usepackage{mathtools}
\DeclarePairedDelimiter{\abs}{\lvert}{\rvert}

\newcommand{\figdir}{eps}

\newcommand{\D}{{D}}
\newcommand{\E}{{\mathcal E}}
\newcommand{\F}{{\mathcal F}}
\newcommand{\Nei}{{\mathcal N}}

\newcommand{\ZZ}{{\mathbb  Z}}
\newcommand{\RR}{{\mathbb  R}}
\newcommand{\zero}{{\textbf 0}}

\newcommand{\dd}{{\bm d}}
\newcommand{\ee}{{\bm e}}
\newcommand{\xx}{{\bm x}}
\newcommand{\yy}{{\bm y}}
\newcommand{\xd}{\bm x'}  
\newcommand{\uu}{{\bm u}}

\newcommand{\SSS}{\mathrm S}
\newcommand{\HHH}{\mathrm H}
\newcommand{\TTT}{\mathrm T}
\newcommand{\OOO}{\mathrm O}

\newcommand{\OO}{\mathrm O}
\newcommand{\PP}{\mathrm P}
\newcommand{\QQ}{\mathrm Q}
\newcommand{\R}{\mathrm R}

\newcommand{\DD}{\Delta}
\newcommand{\EE}{\E}
\newcommand{\s}[2]{{#1}_{#2}}
\newcommand{\DDs}[1]{\s{\DD}{#1}}
\newcommand{\EEs}[1]{\s{\EE}{#1}}

\title{
    Projected images of the Sierpinski tetrahedron and other layered fractal imaginary cubes 
}

\author{Hideki Tsuiki}
\date{Graduate School of Human and Environmental Studies\\ Kyoto University\\
tsuiki@i.h.kyoto-u.ac.jp}

\begin{document}
\maketitle
\begin{abstract}

The Sierpinski tetrahedron has a remarkable property: It is projected to squares in three orthogonal directions, and moreover, to sets with positive Lebesgue measures in numerous directions.
This paper proposes a method for 
characterizing directions along which 
the Sierpinski tetrahedron and other 
similar fractal 3D objects are projected to sets with positive measures.
We apply this methodology to layered fractal imaginary cubes and achieve a comprehensive characterization for them. 
Layered fractal imaginary cubes
are defined as attractors of iterated function systems with layered structures, and they
are  projected to squares in three orthogonal directions. Within this class, the Sierpinski tetrahedron, T-fractal, and H-fractal stand out as exemplary cases.

  \end{abstract}

\section{Introduction}

\begin{figure}[t]
\centering
\subfloat[]
{\includegraphics[height=4.5cm]{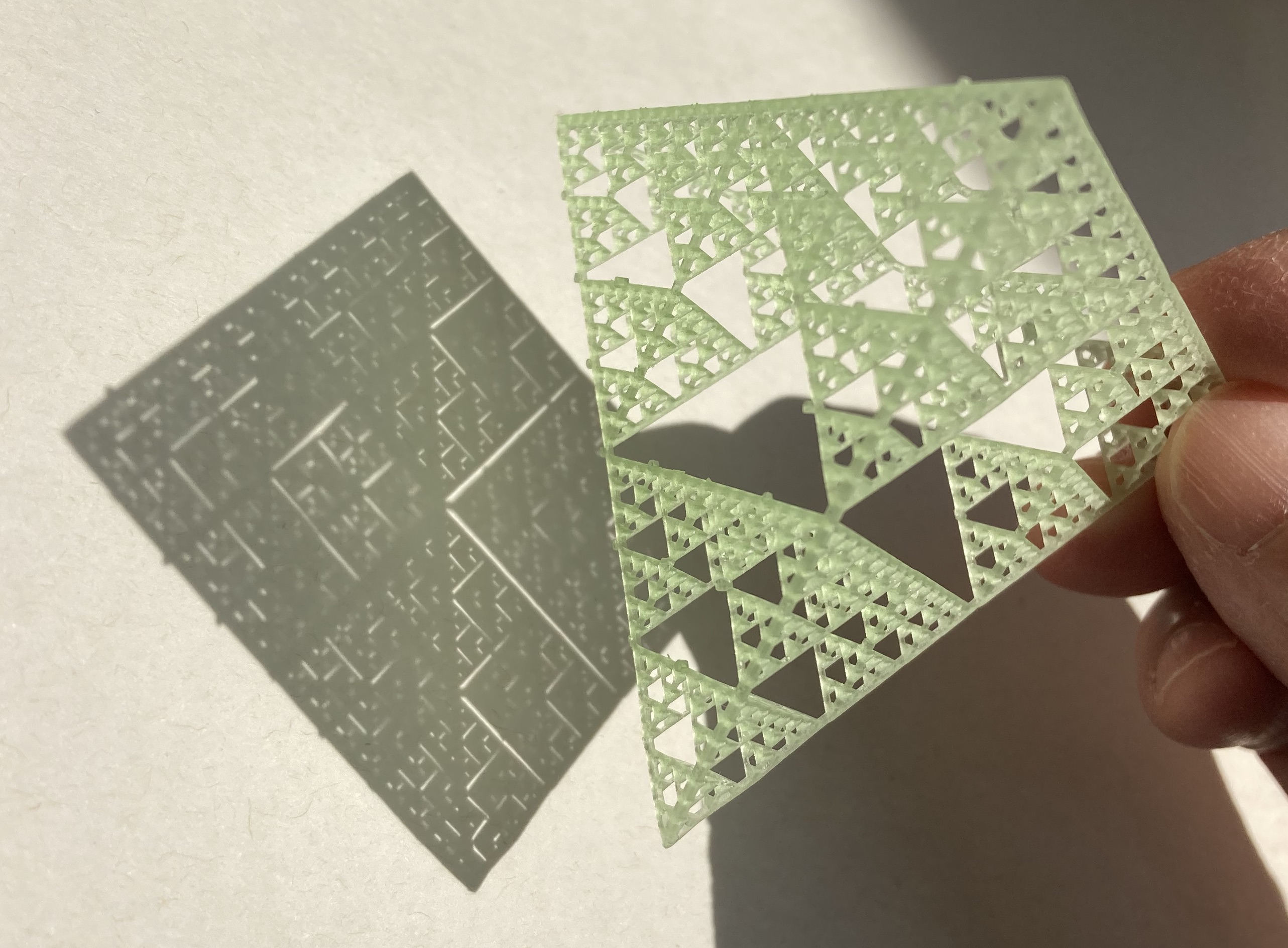} }
\qquad
\subfloat[]{\includegraphics[height=4.5cm]{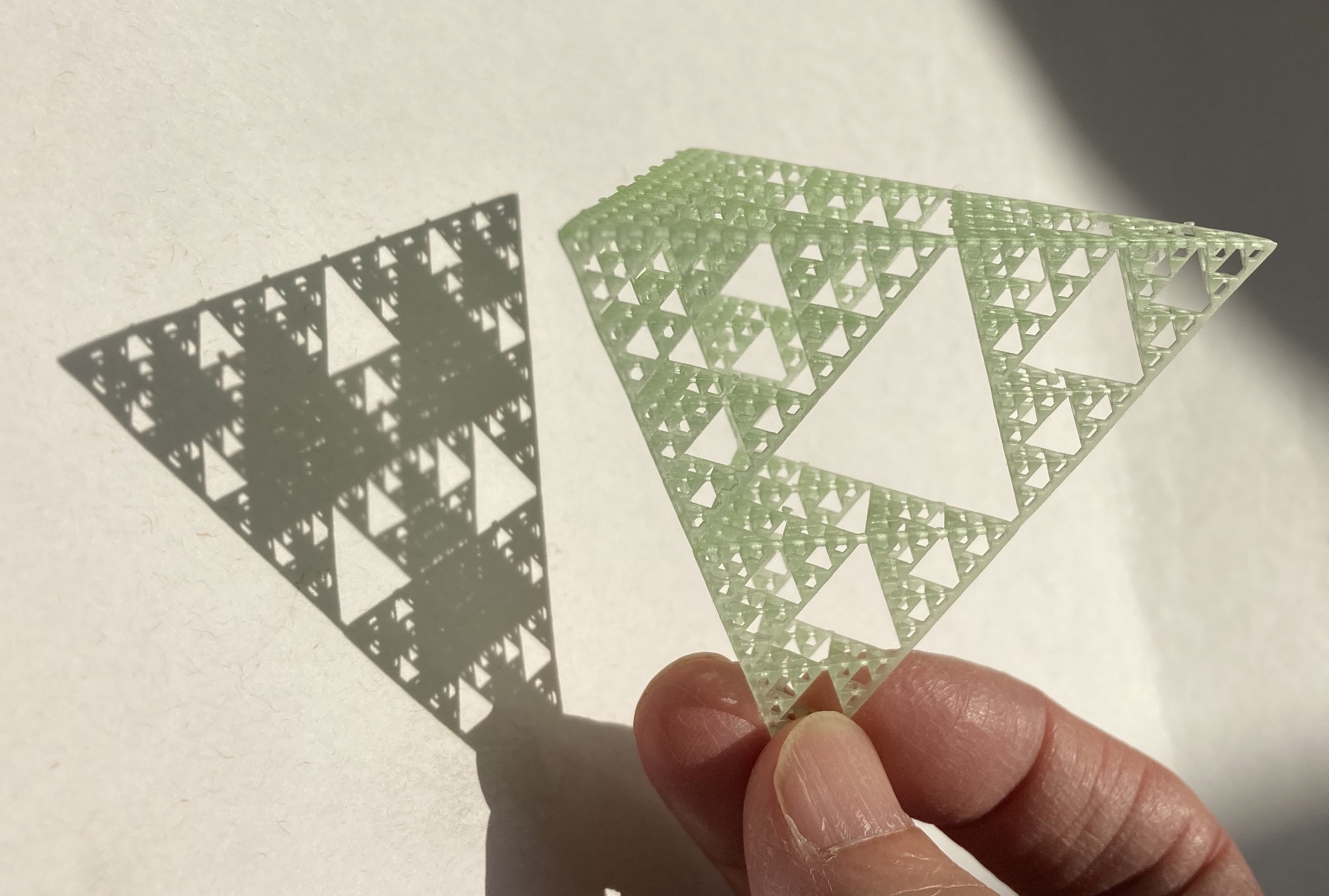}}
\caption{\label{p1}Projected images of the Sierpinski tetrahedron.} 

\end{figure}

The Sierpinski tetrahedron (i.e., three-dimensional Sierpinski gasket) is a three-dimensional fractal object defined as the attractor of 
the iterated function system (IFS) $\{f_\dd(\xx) = 
\frac{\xx + \dd}{2}  \mid \dd \in D_\SSS\}$ where
$D_\SSS \subset \RR^3$ is the set of four vertices of a regular tetrahedron.
That is, it is the unique non-empty compact set $\SSS_\infty$ satisfying 
$\SSS_\infty = \cup_{\dd \in {D_\SSS}} f_\dd(\SSS_\infty)$ 
that exists due to the theory of self-similar fractals by Hutchinson \cite{hutchinson1981}. See Figure~\ref{fig5}(a) for the first levels of its approximations.

The Sierpinski tetrahedron exhibits a remarkable property: when projected along the three orthogonal directions that connect the midpoints of opposite edges, it forms squares, as illustrated in Figure~\ref{p1}(a). Moreover, as demonstrated in Figure~\ref{p1}(b), projections along numerous other directions result in sets with positive Lebesgue measures that tile the  
plane. This paper aims to characterize the specific directions along which the Sierpinski tetrahedron, as well as other similar fractal 3D objects, are projected to sets with positive measures.

As an initial observation, the IFS defining the Sierpinski tetrahedron does not incorporate rotational transformations, and as a consequence, its projected images are also fractals generated by similar IFSs. More precisely, 
we denote by $\F^n(k, D)$ the fractal generated 
by the IFS $\{f_\dd(\xx) = \frac{\xx + \dd}{k} \mid \dd \in \D\}$
for a set of points $D \subset \RR^n$ called the digit set. %
Then, the Sierpinski tetrahedron is denoted by $\F^3(2, D_\SSS)$ and
its image by a projection $\varphi$ is represented as 
$\F^2(2, \varphi(D_\SSS))$.
In general, one can see that a 3D fractal $\F^3(k, D)$ is projected by $\varphi$ to a
2D fractal $\F^2(k, \varphi(D))$.
Therefore, we restrict our attention to fractal objects of the form 
$\F^3(k, D)$ so that we can apply the theory of fractals on their projected images.

Among such fractals, we focus on cases where the digit set is defined as
\[
D_{k, l} = \{(x, y, z)\in \ZZ^3 \mid 0 \leq x, y, z \leq k-1, \ x + y + z \equiv l - 1 \bmod k\}
\]
for integers 
$k \geq 2$ and $0 \leq l < \frac{k}{2}$.
We call the induced fractal object $\F^3(k, D_{k, l})$
a layered fractal imaginary cube of degree $k$.
As we are going to see, it is an imaginary cube 
\cite{algorithm, DBLP:conf/jcdcg/TsuikiT13}, that is, 
it is projected to squares along three orthogonal directions just as a cube.
We call it layered because the digit set is arranged on parallel planes.
$D_{2,0}$ is the set of vertices of a regular tetrahedron. Therefore,
we can choose $\D_\SSS = D_{2,0}$ and the Sierpinski tetrahedron  is the only layered fractal imaginary cube of degree 2.
There are two layered fractal imaginary cubes of degree 3:  
$\F^3(3, D_{3,0})$ that we call the T-fractal,  and $\F^3(3, D_{3,1})$ that we call the H-fractal
(\cite{bridges1,   DBLP:conf/jcdcg/TsuikiT13}, see Figure~\ref{fig5}(b, c)).
The convex hull of the T-fractal is the octahedron T in Figure~\ref{fig2}(b),
which is the convex hull of $\frac{D_{3,0}}{2}$.
Similarly,  
the convex hull of the H-fractal is the dodecahedron H in Figure~\ref{fig2}(c),
which is the convex hull of $\frac{D_{3,1}}{2}$.

The H-fractal is distinguished by its six-fold symmetry, whereas other layered fractal imaginary cubes exhibit three-fold symmetries.
Our main theorem (Theorem~\ref{mainth1}) reflects this difference of symmetries: Except for the H-fractal,
a layered fractal imaginary cube of degree $k$ is projected to a set with positive measure if and only if the projection is
done along $(a, b, c)$ for coprime integers $a, b, c$ such that $a+b+c$ is coprime to  $k$,
and the H-fractal has additional directions that correspond to extra symmetries.

We study measures of projected images through the
theory of self-affine tiles.
Positively measured fractals of the form $\F^2(k, D)$ with  $\abs{D} = k^2$ are
special cases of self-affine tiles,
which have been the subject of
extensive study since the 1990s \cite{Bandt1991, BandtGraf,
  Kenyon1992, Vince1993, Vince1995,  LagariasWang, 
Chun-Kit2017,Fu2015,Su2022,Deng2023
}.
For recent developments in self-affine tiles, refer to \cite{TilingBook}. 
Kenyon demonstrated in \cite{Kenyon1992} that for a four-point set $\D = \{\mathrm{O, P, Q, R}\}$,
the fractal $\F^2(2, D)$ has 
positive Lebesgue measure if and only if $ r\vec{\mathrm{OR}} = p \vec{\mathrm{OP}} + q
\vec{\mathrm{OQ}}$ for odd integers $p, q, r$.
This result leads to a characterization of  
positively measured 
projections of the Sierpinski tetrahedron,
which is consistent with
the characterization provided by Theorem \ref{mainth1}. %
By extending Kenyon's proof concepts and utilizing the notion of projection of 
differenced radix expansion sets, we
present
a method for characterizing directions along which fractal objects are projected to sets with positive Lebesgue measures.
Then, we apply it to layered fractal imaginary cubes.
    
Kenyon also studied one-dimensional projections of the one-dimensional Sierpinski Gasket in \cite{Kenyon1997}, characterizing the directions along which it is projected to  sets with positive Lebesgue measures. Our work provides similar characterization for two-dimensional projections of a series of three-dimensional fractal objects.

In the next section, 
we introduce the concept of layered fractal imaginary cubes 
and explain our results in Section \ref{sec:mainresults}.
Section \ref{affine} provides an overview of the properties of self-affine tiles
relevant to our analysis and propose our method.  Using this methodology, 
the results in Section \ref{sec:mainresults}
are proved in Section \ref{section:proof}.
In section \ref{sec:conclusion},
we discuss the case of projections of
non-layered imaginary cubes.  In Section \ref{sec:conclusion},
we also explore fractal imaginary `squares', and explain its relation to
Kenyon's work on 
projections of the one-dimensional Sierpinski gasket.

\section{Layered fractal imaginary cubes}\label{section-imaginarycube}

\begin{figure}[t]
  \centering
  \subfloat[{{Regular tetrahedron}}]
  {\qquad\includegraphics[width=25mm]{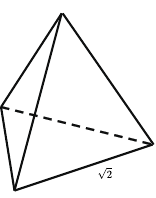}\qquad\qquad}
  \subfloat[T]
  {\includegraphics[width=25mm]{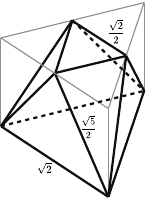}}\qquad\qquad
  \subfloat[H]
  {\includegraphics[width=27mm]{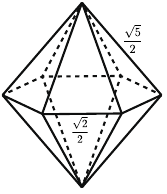}}
  \caption{Examples of imaginary cubes. \label{fig1}} 
  \bigskip
  
\centering
\subfloat[{{Regular tetrahedron}}]
{\qquad\includegraphics[width=27mm]{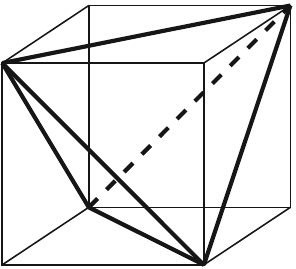}\qquad\qquad}
\subfloat[T]
{\includegraphics[width=27mm]{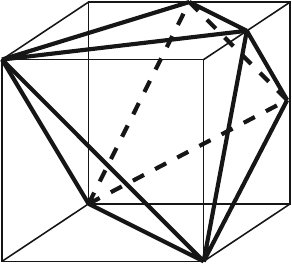}}
\qquad\qquad
\subfloat[H]
{\includegraphics[width=27mm]{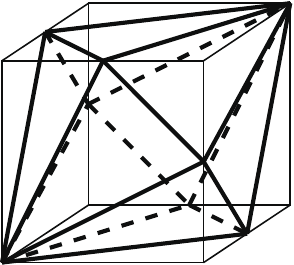}}\qquad

\caption{Objects in Figure~\ref{fig1} located in unit cubes.\label{fig2}} 
\end{figure}

We start by %
introducing 
the notion of imaginary cubes,
introduced and developed in %
\cite{bridges1, algorithm, DBLP:conf/jcdcg/TsuikiT13}.

\begin{definition} \label{Def1}
    An {\em imaginary cube} is an object that is projected to squares in  three orthogonal directions, and in each of these projected squares, the four edges are parallel to the other two projection directions.
\end{definition}
\noindent 
Note that a regular octahedron also has square projections in three orthogonal directions, but
it is not an imaginary cube.
An imaginary cube is contained in a cube $C$ defined by the three square projections.
When it is necessary to specify $C$,  we say that it is an imaginary cube of a cube $C$.

A regular tetrahedron is an example of an imaginary cube.
There are two more polyhedral imaginary cubes that are relevant to this paper: 
triangular antiprismoid imaginary cube (abbreviated as T)
and hexagonal bipyramid imaginary cube (abbreviated as H). 
 They are defined in  Figure~\ref{fig1}, and cubes of which they are imaginary cubes are illustrated in Figure~\ref{fig2}.
T and H
exhibit remarkable properties:
H is a {\em double imaginary cube} in that it is an imaginary cube of two different cubes and thus it has six square projections.  
T is notable in that 
its three diagonals are orthogonal to each other.
Furthermore,  
T and H
 together form a three-dimensional tiling of the space as explained in \cite{algorithm}.

These three imaginary cubes -- the regular tetrahedron, H and T -- share the common
property of being convex hulls of fractal imaginary cubes.
According to Hutchinson \cite{hutchinson1981},
for contractions $f_i: \mathbb{R}^n \to \mathbb{R}^n$ $(i = 1,\dotsc,m)$, 
an IFS (iterated function system) $I = \{f_i \mid i = 1,2,\ldots,m\}$ 
defines a self-similar fractal object $\mathcal{F}_I$
as the fixed point of 
the following contraction map on
the metric space ${\mathcal H}^n$ of non-empty compact subsets of $\mathbb{R}^{n}$.
\begin{align*}
F(X) &= \bigcup_{i=1}^{m} f_{i}(X).%
\end{align*}
For any $A \in {\mathcal H}^n$, the sequence 
$(A_i)_{i\geq 0}$ where $A_0 = A$ and 
$A_{i+1} = F(A_i)$
for $i\geq 0$ 
converges to $\mathcal{F}_I$ and $F(\mathcal{F}_I) = \mathcal{F}_I$.

In this article, we only consider the case
an IFS takes the form 
\[\{f_\dd(x) = \frac{x+\dd}{k} \mid \dd \in \D\}\]
for some 
integer $k \geq 2$ and some set $\D \subset \RR^n$ called the digit set.  That is, the components are similarity maps that do not perform rotational transformations. 
In this case, the map $F$ can be expressed as
\begin{align*}
F(X) &= \frac{X + D}{k} %
\end{align*}
where $+$ denotes the Minkowski sum
$X + Y = \{\xx + \yy \mid \xx \in X, \yy \in Y\}$
 and $\frac{X}{k}$ denotes $\{\frac{\xx}{k} \mid \xx \in X\}$. 
We denote the induced fractal object as 
$\F^n(k, D)$ or $\F(k, D)$ if the dimension is not important.

With this restriction, the image of the fractal
$\F^3(k, D)$ by a projection $\varphi$
is a two-dimensional fractal $\F^2(k, \varphi(D))$.
To ensure that $\F^2(k, \varphi(D))$ has a positive Lebesgue measure, it is necessary that its Hausdorff dimension is equal to 2.
On the other hand, Hausdorff dimension cannot exceed the similarity dimension
$\log_k \abs{D}$ of the iterated function system. 
In this paper, we focus on the case where $|D| = k^2$, representing fractals with the minimum number of IFS components. We call $\F^3(k, D)$ under the condition $|D| = k^2$ a 'homothetic %
fractal imaginary cube' of degree $k$, or simply a 'fractal imaginary cubes' of degree $k$, because all 
the fractal and imaginary cube objects discussed in this paper conform to this form.

\begin{lemma}[\cite{bridges1}]\label{lemmal1}
  Let $D$ be a three-dimensional digit set of cardinality $k^2$.
  $\F^3(k, \D)$ is a fractal imaginary cube of degree $k$
  if and only if, for some cube $C$, 
  $F(C)$ for $F(X) = \frac{X+D}{k}$  
  is an imaginary cube of $C$.
\end{lemma}
\begin{proof}    
If $\F^3(k, \D)$ is an imaginary cube, then there is a cube $C$ of which
$\F^3(k, \D)$ is an imaginary cube.  
Since $f_\dd(\xx) = \frac{\xx+\dd}{k}$ fixes $\frac{\dd}{k-1}$,
$\frac{\dd}{k-1}$ for $\dd \in D$ are contained in $\F^3(k, \D) \subseteq C$.
Since $f_\dd(\xx)$ is on the line segment between $\xx$ and $\frac{\dd}{k-1}$ and $C$ is convex, 
$F(C) \subseteq C$. Since, in addition, $X \supseteq Y$ implies $F(X) \supseteq F(Y)$,
The sequence $(A_i)_{i \geq 0}$ for $A_0 = C$ and $A_{i+1} = F(A_i)$
forms a decreasing sequence
$C = A_0 \supset A_1 \supset A_2 %
\supset
\ldots$ whose intersection is $\F^3(k, \D)$.
Therefore, $A_i$ are all imaginary cubes of $C$ because
$C$ and $\F^3(k, \D)$ have the same three square projections.
In particular, $A_1 = F(C)$ is an imaginary cube of $C$. 

Conversely, for a cube $C$, suppose that $F(C)$ is an imaginary cube of $C$. 
Let $(A_i)_{i \geq 0}$ be the sequence
defined as $A_0 = C$ and $A_{i+1} = F(A_i)$ and $\varphi$ be a projection %
to a square $S$ 
along an edge of $C$. 
We show by induction that $A_{i}$ 
($i = 0, 1, \ldots$)
are all projected to $S$ by $\varphi$.
First, we have $\varphi(A_0) = S$.
Suppose that $\varphi(A_i) = S$.
Then, since $\varphi(A_i) = \varphi(C)$, we have
$\varphi(f_\dd(A_{i})) = \varphi(f_\dd(C))$. Therefore, 
both $A_{i+1} = \cup_{\dd \in D} f_\dd(A_{i})$ and
$A_{1} = \cup_{\dd \in D} f_\dd(C)$ are projected by $\varphi$ to the same set, which is $S$
because $A_1$ is an imaginary cube of $C$. 
Thus, $A_{i}$ ($i = 0, 1, \ldots$) are all projected by $\varphi$ to $S$. 
Then, $\F^3(k, \D)$ is also projected to 
$S$
because projection $\varphi$ is a continuous map from $\mathcal{H}^3$ to 
$\mathcal{H}^2$.
\end{proof}

If $F(C)$ is an imaginary cube, then 
it is a union of $k^2$ cubes
selected from the $k^3$ cubes obtained by cutting $C$ into $k \times k \times k$ small cubes so that they do not overlap when viewed from the three face-directions of $C$.
Such a selection exists corresponding to a Latin square of size $k$,
which is a $k \times k$ matrix of $\{0, 1,2,\ldots,k-1\}$ 
with each number appearing exactly once in each row and column.

\begin{figure}[t] 
\centering
\includegraphics[width=12cm]{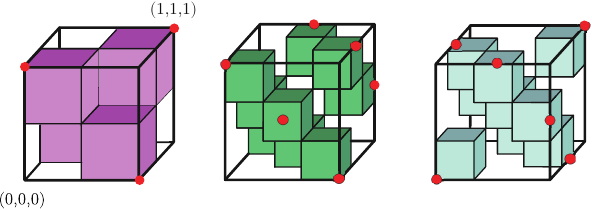}\\ \vspace*{-0.2cm}
\hspace*{0.2cm}({\textbf a}) $D_{2,0}$   \hspace*{2.5cm}       ({\textbf b}) $D_{3,0}$    \hspace*{3cm}  ({\textbf c}) $D_{3,1}$
\hspace*{0.5cm}
\caption{\label{figfrac}
The arrangements of cubes by (a) $\D_{2, 0}$, 
(b) $\D_{3, 0}$, and (c) $\D_{3, 1}$, which
generate $\SSS_\infty$, $\TTT_\infty$, and $\HHH_\infty$, respectively. %
Points in $\frac{D_{k,l}}{k-1}$ are marked with red circles.
}
\end{figure}

When $k = 2$, there are two such arrangements of $4$ cubes but
they are congruent (Figure~\ref{figfrac}(a)),
and therefore there is only one fractal imaginary cube of degree 2. 
 When $k=3$, there are two arrangements of $9$ cubes modulo congruence as depicted in Figure~\ref{figfrac}(b) and (c),
 and therefore there exist two fractal imaginary cubes of degree 3, which 
  we will investigate later.
For 
$k \geq 4$, there exist 36 fractal imaginary cubes of degree 4, and 3482 fractal imaginary cubes of degree 5, modulo congruence \cite{bridges1}. 
Since the digit sets of fractal imaginary
cubes of degree $k$ correspond to Latin squares of size $k$ and it is demonstrated in \cite{lint_wilson_2001}
that the number of Latin squares of size $k$ is lower bounded by $\frac{(k!)^{2k}}{k^{k^2}}$,
the number
of fractal imaginary cubes increases rapidly with $k$.

Among them, there is a class of digit sets
with a uniform structure.
For integers %
$k \geq 2$ and $0 \leq l < \frac{k}{2}$,
we define the digit set $D_{k, l}$ as 
\[
D_{k, l} = \{(x, y, z)\in \ZZ^3 \mid 0 \leq x, y, z \leq k-1, \ x + y + z \equiv l - 1 \bmod k\}.\label{e:layered}
\]
For the unit cube $C$, $\frac{C + D_{k,l}}{k}$ is a union of $k^2$ cubes that form an imaginary cube. Therefore,
by Lemma \ref{lemmal1},  
$\F(k,D_{k, l})$ is a fractal imaginary cube, which we call a  {\em layered fractal imaginary cube}.  
Though $D_{k,l}$ can be defined for $0 \leq l < k$, $D_{k, l}$ and $D_{k, k-l-1}$ are congruent and thus we only consider
$0 \leq l < \frac{k}{2}$.
The set $D_{k, 0}$ is contained in the union of two planes
$x + y + z = -1$ and 
$x + y + z = k -1$, whereas
$D_{k, l}$ for $0 < l < \frac{k}{2}$ is contained in the union of three planes:
$x + y + z = l-1$,
$x + y + z = k +l-1$ and 
$x + y + z = 2k + l-1$. 
In this paper, we investigate the projections of layered fractal imaginary cubes
and conclude with a discussion on projections of fractal imaginary cubes in the final section.

In a fractal imaginary cube $\F(k,D)$, the set of fixed points of components of the IFS is $\frac{D}{k-1}$, which is 
necessarily contained in $\F(k,D)$.
On the other hand, for $P$ the 
polyhedron obtained as the
convex hull of $\frac{D}{k-1}$, $P \supset \F(k,D)$
because \[P \supset F(P)
\supset F(F((P))) \supset \cdots\] forms a decreasing sequence whose intersection is 
$\F(k,D)$, as we proved for
a cube $C$ containing $\F(k,D)$ in Lemma \ref{lemmal1}. Thus, we have the following:
\begin{lemma}
  The convex hull of $\F(k,D)$ is the polyhedron obtained as the convex hull of $\frac{D}{k-1}$.
\end{lemma}

In Figure~\ref{figfrac}(a), the points of $D_{2, 0}$ are marked with red circles, 
whose convex hull is a regular tetrahedron.  Since it is the convex hull of the fractal
$\F(2, D_{2, 0})$, $\F(2, D_{2, 0})$ is the Sierpinski tetrahedron
which we denote by $\SSS_\infty$. 

\begin{figure}[t] 
\centering
\begin{tabular}[b]{cc}
\raisebox{1cm}{ (\textbf a)} &\includegraphics*[height=30mm]{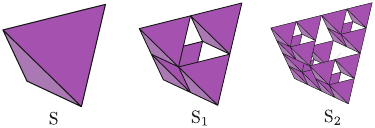}\\
\raisebox{1cm}{(\textbf b)}&\includegraphics*[height=30mm]{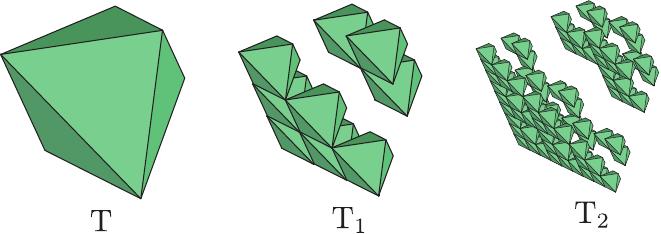}\\
\raisebox{1cm}{(\textbf c)} &\includegraphics*[height=30mm]{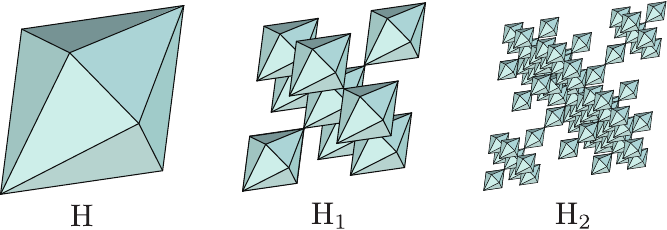}
\end{tabular}
\caption{\label{fig5}
The first two approximations of 
(a) $\SSS_\infty$, (b) $\TTT_\infty$, and (c) $\HHH_\infty$.}
\end{figure}

For $D_{3, 0}$ and $D_{3, 1}$, the sets of the fixed points 
of the IFS components are 
$\frac{D_{3, 0}}{2}$ and $\frac{D_{3, 1}}{2}$, respectively, as indicated by red circles in Figure~\ref{figfrac}(b) and (c).
The convex hulls of these sets are
imaginary cubes T and H, respectively.
Therefore, $\F(3, D_{3,0})$ and $\F(3, D_{3,1})$ are 
fractal imaginary cubes 
with their convex hulls being T and H, respectively.
These fractals are denoted by $\TTT_\infty$ and
$\HHH_\infty$, respectively.
Figure~\ref{fig5} illustrates the first two iterations  of approximations of $\SSS_\infty$, $\TTT_\infty$, $\HHH_\infty$
starting from their convex hulls.
We define the shapes of $\TTT_\infty$ and $\HHH_\infty$ as T-fractal and H-fractal, respectively.

Note that only the H-fractal exhibits six-fold symmetry and is a
double imaginary cube, whereas the other layered fractal imaginary cubes exhibit three-fold symmetry along the vector $(1,1,1)$.
This extra symmetry of the H-fractal leads to different type of characterization
of projections that yields images with positive measure in Theorem~\ref{mainth1}.

\section{Main results}\label{sec:mainresults}

In this section, we state our results about projected images of
layered fractal imaginary cubes.
As noted,  
projected images of layered fractal imaginary cubes also possess fractal structures.
The appendix contains photographs of 3D printed models of three representative 
layered fractal imaginary cubes; Sierpinski tetrahedron, T-fractal, and H-fractal,  
casting shadows under sunlight
from typical directions.
We have marked with red frames those photographs that show shadows with positive measures even when the
models were perfect mathematical representations.

\begin{figure}[t]
\centering
\subfloat[$\SSS_\infty$, $(1, 1, 1)$]
{\includegraphics*[height=25mm]{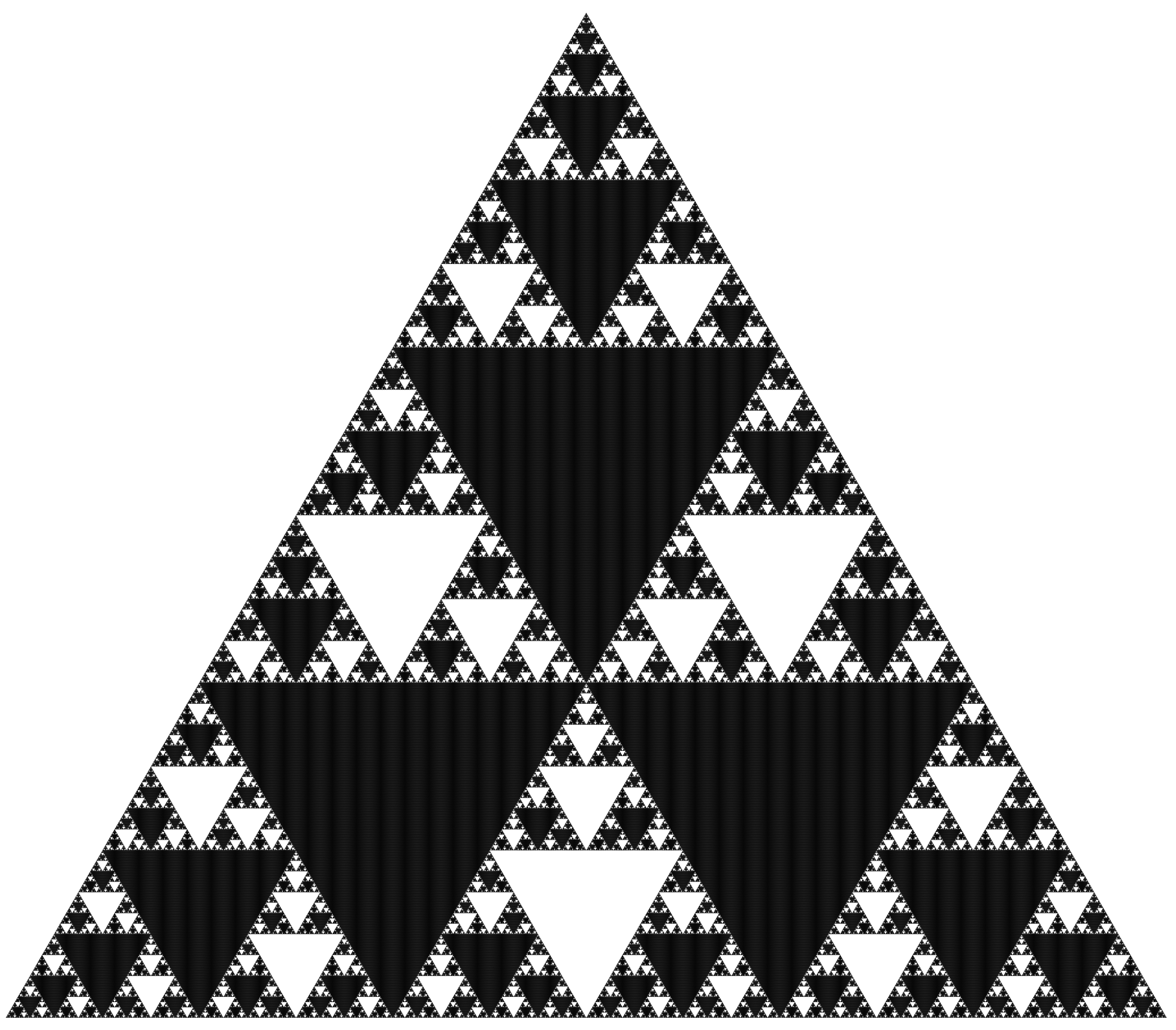}}\quad
\subfloat[$\TTT_\infty$, $(1, 1, -1)$ 
]
{\includegraphics*[height=25mm]{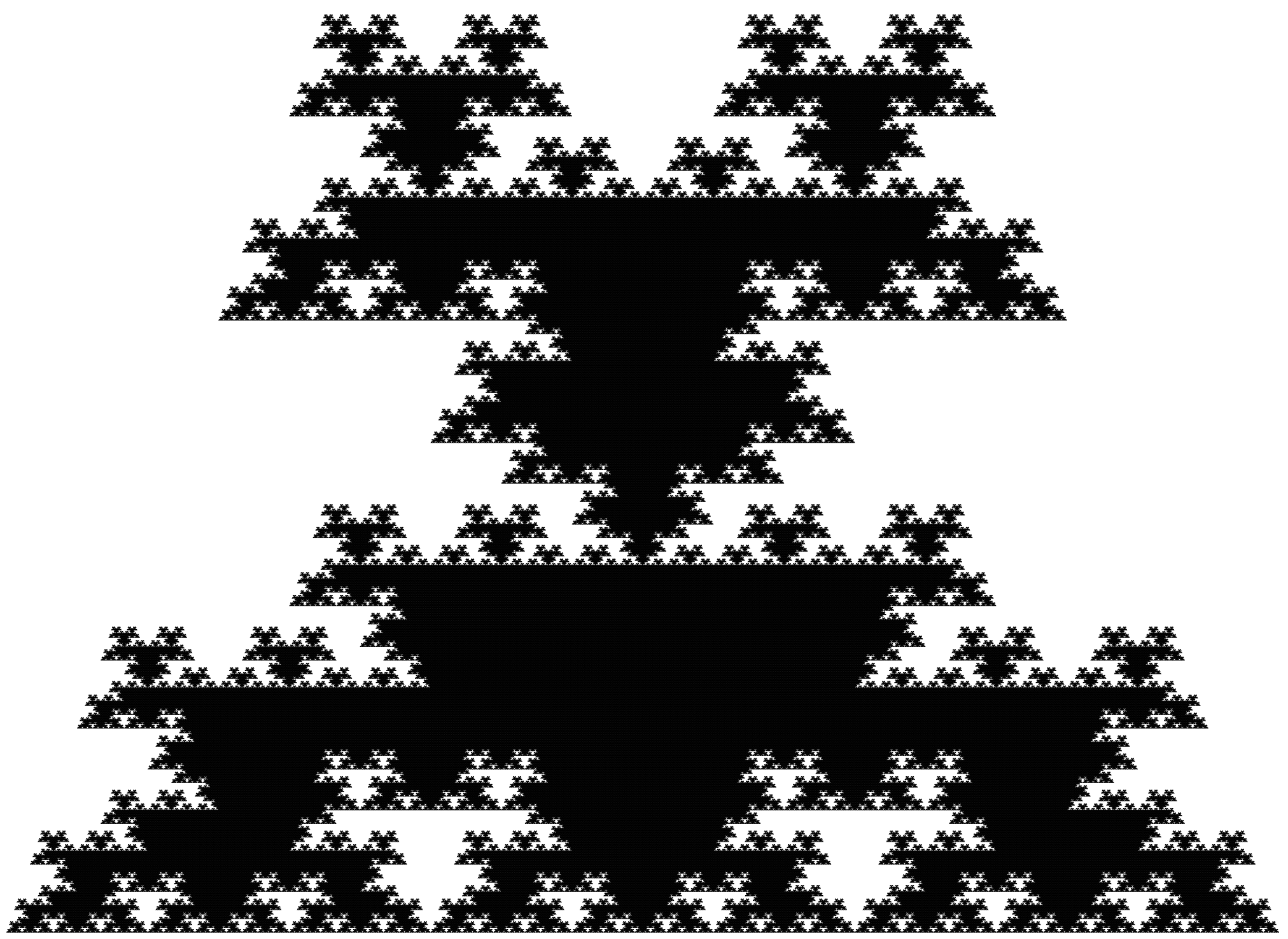}}\quad
\subfloat[$\HHH_\infty$, $(1,1,-1)$ %
and $(1,1,-5)$
]
{\includegraphics*[height=30mm]{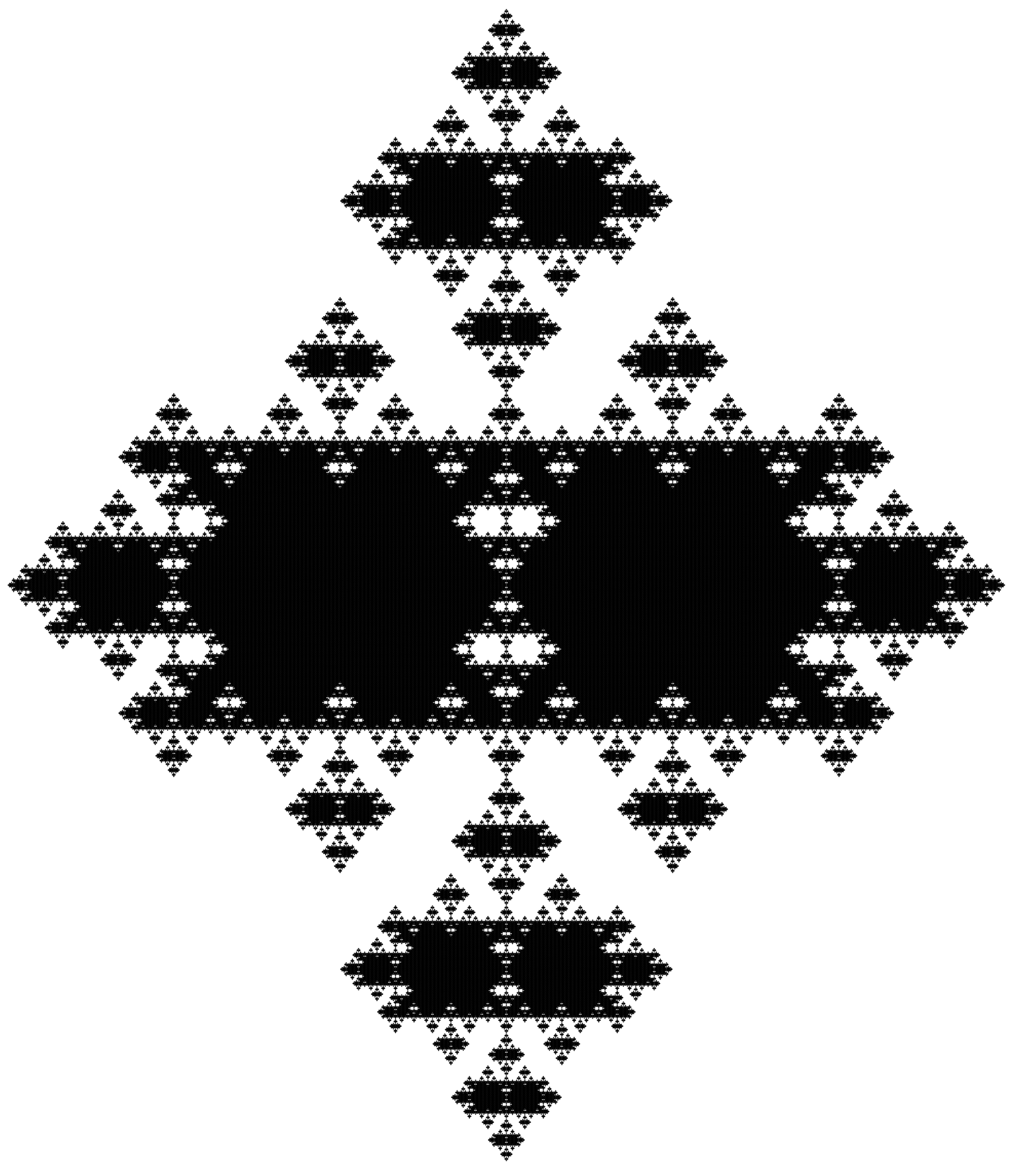}}\\
\subfloat[$\SSS_\infty$, $(1,2,0)$]
{\includegraphics*[height=25mm]{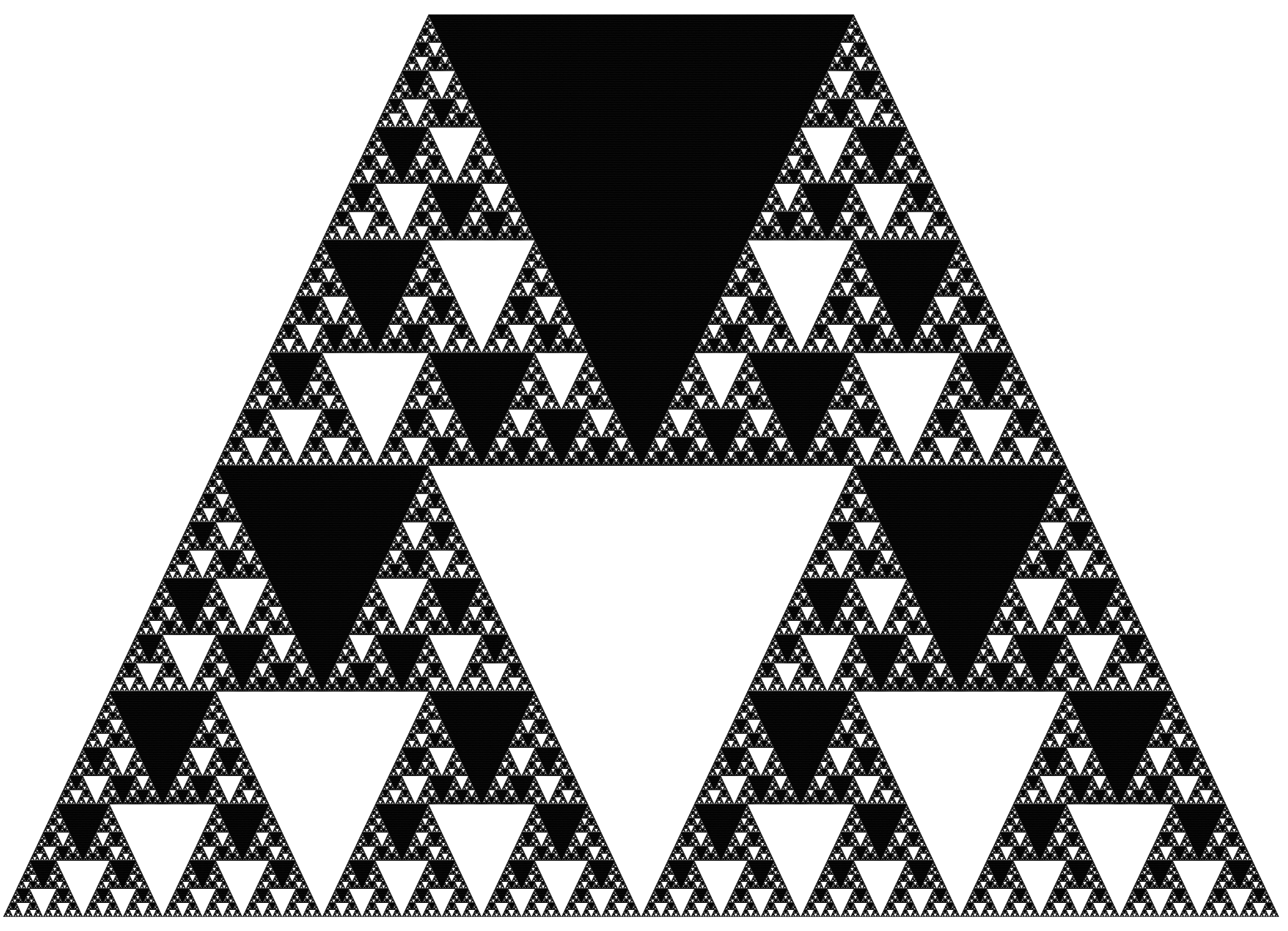}}\quad
\subfloat[$\TTT_\infty$, $(1, 1, 0)$]
{\includegraphics*[height=25mm]{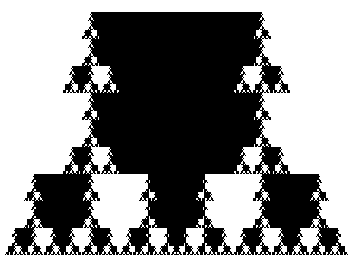}}\quad
\subfloat[$\HHH_\infty$, $(1, 1, 0)$ %
and $(1,1,4)$
]
{\includegraphics*[height=25mm]{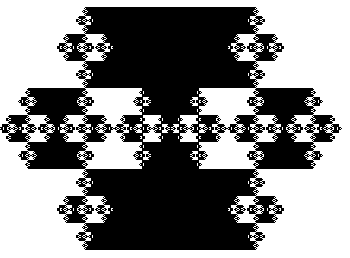}} \\
\subfloat[$\F(4, \D_{4,0})$, $(1,1,1)$]
{\includegraphics*[height=28mm]{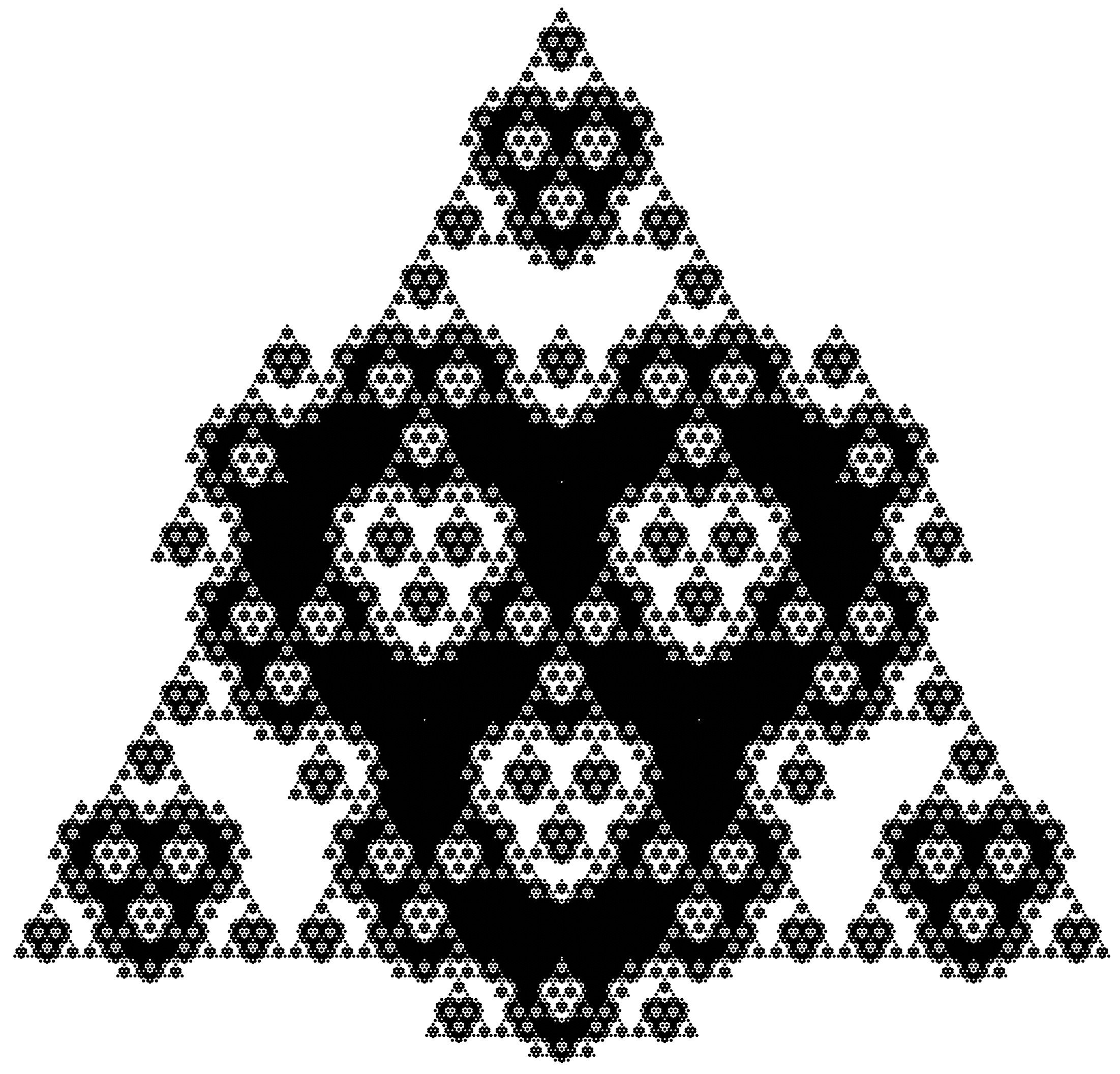}}\quad\quad
\subfloat[$\F(4, \D_{4,1})$, $(1,1,1)$]
{\includegraphics*[height=28mm]{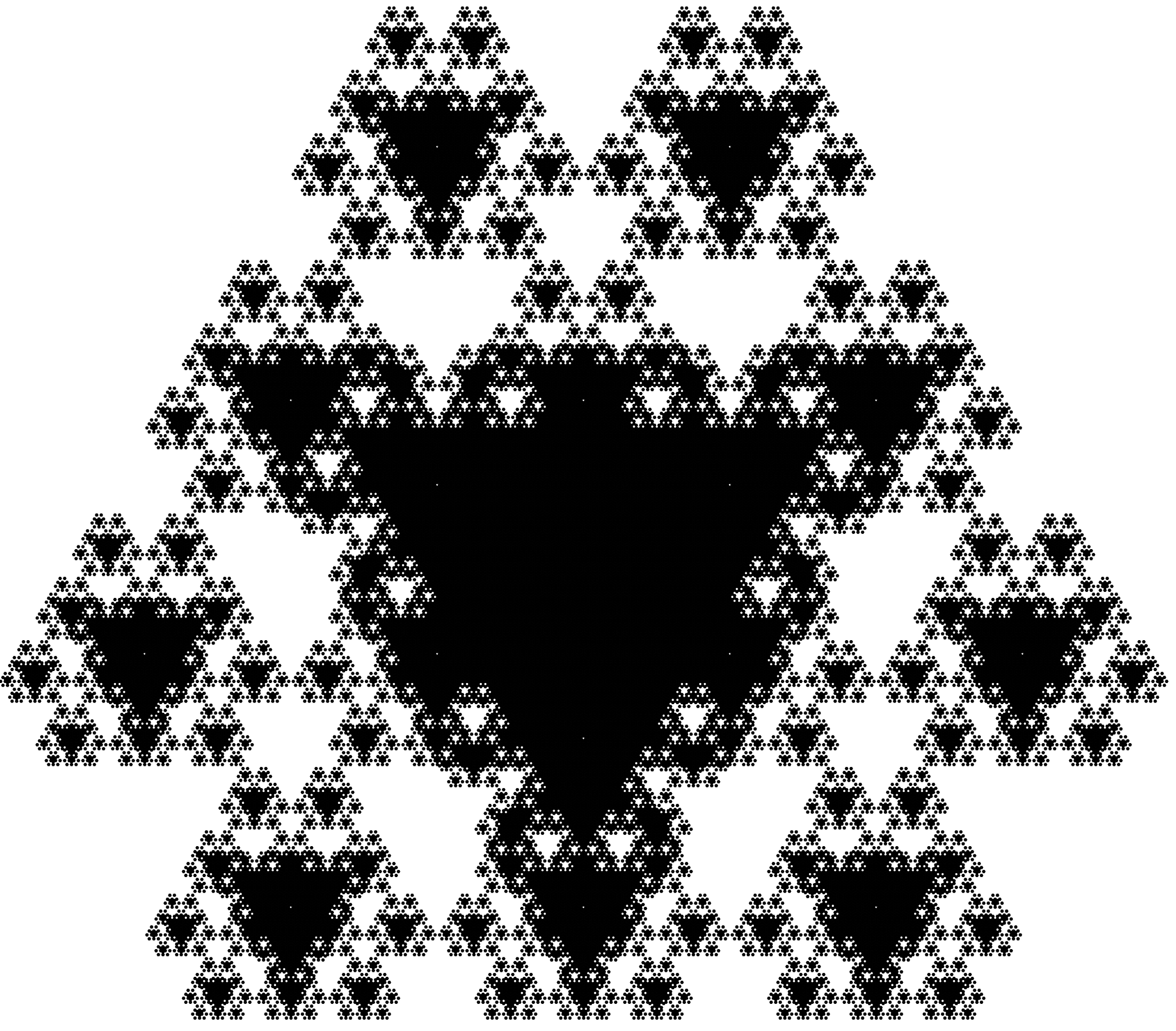}}\quad\quad
\subfloat[$\F(5, \D_{5,2})$, $(1,1,1)$]
{\includegraphics*[height=29mm]{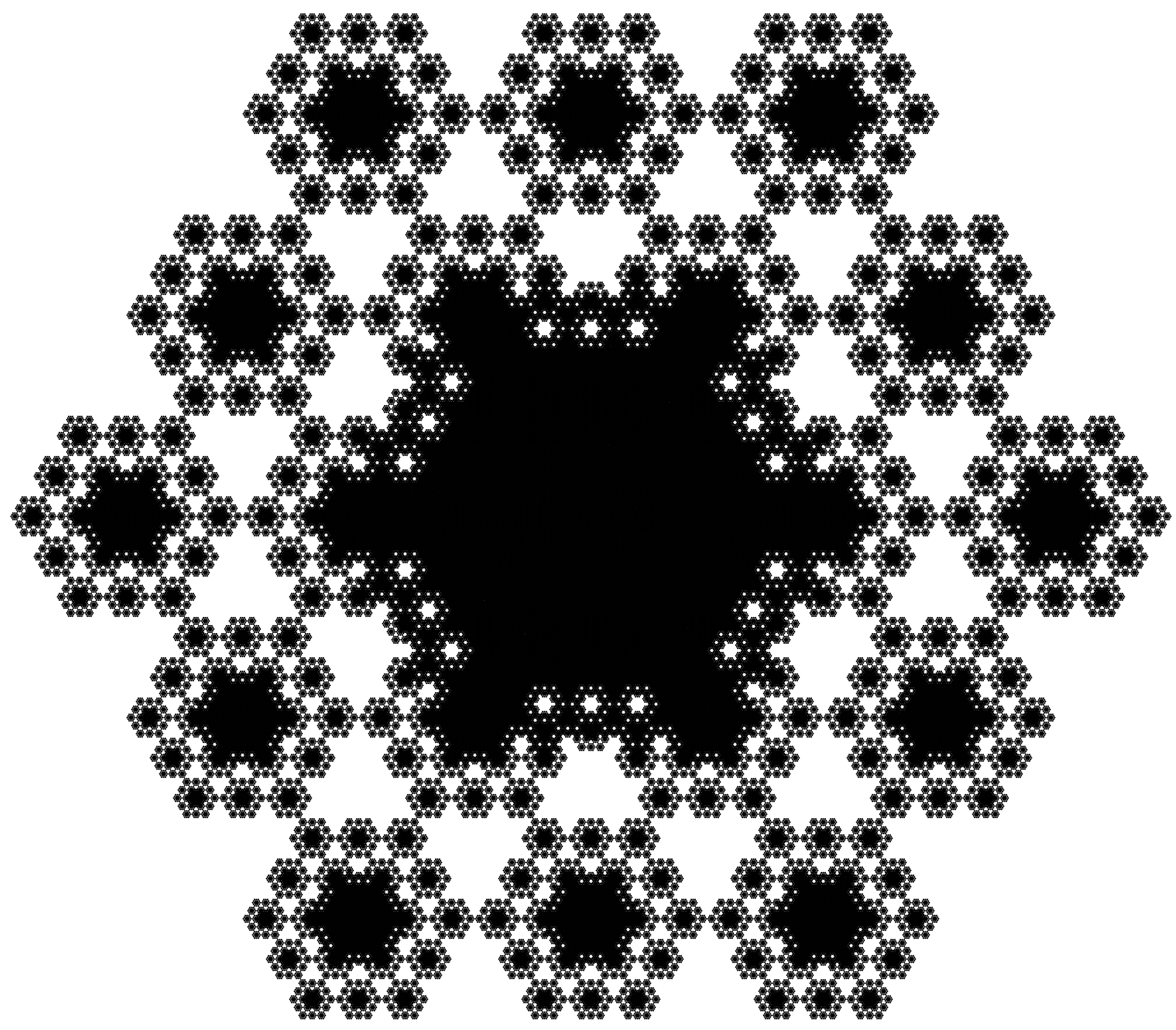}} 
\caption{\label{fig7} Projected images of layered fractal imaginary cubes.
  } 
\end{figure}

Some of the positively measured images of layered fractal imaginary cubes are 
depicted in Figure~\ref{fig7}.
In addition to these images,
all layered fractal imaginary cubes %
have  square projections along the three orthogonal vectors
$(1, 0, 0),
(0, 1, 0)$ and $(0, 0, 1)$ as imaginary cubes. 
Furthermore, $\HHH_\infty$, being a double imaginary cube, is projected to squares  also along vectors $(-1, 2, 2), (2, -1, 2)$, and $(2, 2, -1)$, obtained by rotating the above three vectors 180 degrees around the axis $(1, 1, 1)$.

We now present our main theorem, which characterizes the directions along which layered fractal imaginary cubes are cast into images with positive Lebesgue measures.
We call a vector $(a,b,c) \in \RR^3$ such that $a, b, c$ are coprime integers a {\em coprime vector}.

\begin{theorem}\label{mainth1}
    Let $k \geq 2$ and $0 \leq l < \frac{k}{2}$.
     $\F(k,D_{k,l})$ is projected to a set with positive Lebesgue measure if and only if
     the projection is done along a coprime vector $(a,  b, c)$  such that
   \begin{itemize}
    \item[(a)] $a +  b + c$ is coprime to $k$ \ (when $(k, l) \ne (3,1)$); 
    \item[(b)] $a +  b + c$ is coprime to $3$ or   $a \equiv  b \equiv -\frac{a+b+c}{3} \not \equiv 0 \ (\bmod\ 3) $\ \ 
    (when $(k, l) = (3,1)$, i.e., $\F(k, D_{k,l}) = \HHH_\infty $).
   \end{itemize}
   \end{theorem}

Note that
$a \equiv  b \equiv -\frac{a+b+c}{3} \ (\bmod\ 3)$ implies $a \equiv  b \equiv  c\ (\bmod\ 3)$.
Note also that the 'only if' part of this theorem says that
if (1) projection is done along a coprime vector $(a,b,c)$  that does not meet the conditions (a) and (b) 
or (2) any pair among $a, b$ and $c$ forms an irrational ratio, then the projected image is a null set.

According to Theorem \ref{mainth1}(a), $\SSS_\infty$ and $\TTT_\infty$ are projected to sets with positive Lebesgue measures if $a + b + c$ is not a multiple of 2 and 3, respectively.  As a special case of this theorem, we have:

\begin{corollary}
 The projected image of $\F(k,D_{k, l})$ along the vector $(1,1,1)$
 has a positive measure if and only if $k$ is not a multiple of 3.
\end{corollary}

{\em Note1: }
Since all the layered fractal imaginary cubes exhibit three-fold symmetry around the axis $x=y=z$, their projected images along (1,1,1) also inherit this three-fold symmetry as well. Refer to Figure~\ref{fig7}(a), (g), (h), and (i), as well as to the pictures
in the Appendix showing shadows cast along (1, 1, 1). 
Among them, the projected image of $H_\infty$, known as the hexaflake fractal, exhibits six-fold symmetry.
More generally, six-fold symmetry appears in images of $\F(k, D_{k, \frac{k-1}{2}})$ for each odd number $k$.
See Figure~\ref{fig7}(i) for an example with $k=5$. 
This fact is implied by the hexagonal arrangement of the projected image of $D_{k,l}$ along $(1,1,1)$.

{\em Note2: }
This theorem also says that if the projection is parallel to the plane $x + y + z = 0$, then the image of $\F(k,D_{k, l})$ is a null set for every 
$k$ and $l$.  Refer to the picture of the projected images of 
$\SSS_\infty$, $\HHH_\infty$, and
$\TTT_\infty$ along $(1, 1, -2)$ in the Appendix.
One can easily see that if $(k, l)$ is not $(2, 0)$ nor $(3, 1)$, then 
such a projected image is contained in a union of parallel lines. In particular, 
for the case of projection of $\TTT_\infty$ along $(1, 1, -2)$, the image is a union of line segments arranged according to the Cantor set. 

As stated in Theorem \ref{mainth1}(b), the directions along which 
$\HHH_\infty$ is projected to sets with positive measures are divided into two disjoint sets.
The first one is the same as that of $\TTT_\infty$  and the second one is obtained
by rotating the first one by 180 degrees around the axis $(1, 1, 1)$.
We first show this fact.

\begin{lemma}\label{lemma-H}
An integer vector $(a, b, c)$ such that
  $3 \nmid a +  b + c$ is  
  rotated by a 180-degree rotation around the axis $(1, 1, 1)$ to  $\frac{1}{3}(\alpha, \beta, \gamma)$ for 
integers $\alpha, \beta, \gamma$ satisfying
  $\alpha \equiv \beta \equiv - \frac{\alpha + \beta + \gamma}{3} \not \equiv 0 \ (\bmod\ 3)$, and vice versa.

\begin{proof}
  Since the vectors $\frac{1}{3}(-1, 2, 2)$, $\frac{1}{3}(2, -1, 2)$ and $\frac{1}{3}(2, 2, -1)$ are obtained by rotating
the vectors $(1, 0, 0)$, $(0, 1, 0)$ and $(0, 0, 1)$ respectively,
we define 
\[
(\alpha, \beta, \gamma) = a (-1, 2, 2) +  b (2, -1, 2) + c (2, 2, -1)
\]
for $a, b, c \in \RR$.
We aim to show that $a, b, c$ are integers such that $3 \nmid a +  b + c$ if and only if 
$\alpha, \beta, \gamma$ are integers such that
  $\alpha \equiv \beta \equiv - \frac{\alpha + \beta + \gamma}{3} \not \equiv 0 \ (\bmod\ 3)$.

  If $a +  b + c \equiv 1 \ (\bmod\ 3)$, 
  then $\alpha + \beta + \gamma =  3(a + b + c) \equiv 3 \mod 9$
  and thus $\frac{\alpha + \beta + \gamma}{3} \equiv 1 \ (\bmod\ 3)$.
  Since $\alpha = 2(a + b + c) - 3a$, it follows that  $\alpha \equiv -1 \ (\bmod\ 3)$.
  Similarly, $\beta$ also satisfy 
$\beta \equiv -1 \ (\bmod\ 3)$.
  In the case where $a +  b + c \equiv -1 \ (\bmod\ 3)$, 
  we have $\alpha + \beta + \gamma \equiv -3 \mod 9$ and therefore
  $\frac{\alpha + \beta + \gamma}{3} \equiv -1 \ (\bmod\ 3)$ and
  $\alpha \equiv \beta \equiv 1 \ (\bmod\ 3)$.

  Conversely, if $\alpha \equiv \beta \equiv - \frac{\alpha + \beta + \gamma}{3} \not \equiv 0 \ (\bmod\ 3)$, then
  $a +  b + c = \frac{\alpha + \beta + \gamma}{3}  \not \equiv 0 \ (\bmod\ 3)$.
  $a = \frac{1}{9
  }(-\alpha + 2\beta+ 2\gamma) =  
  \frac{1}{3}(2\frac{\alpha + \beta + \gamma}{3}  - \alpha)$ is an integer,
  and similarly for $b$ and $c$.
\end{proof}
\end{lemma}

We prove
Theorem \ref{mainth1} via 
affine-transformed images.
We consider the affine transformation 
\[
\psi_{k,l}(x,y,z) = (
x, y, \frac{x+y+z-l+1}{k}-1).
\]
It translates the digit set $D_{k,l}$ to $D'_{k,l} = \psi_{k,l}(D_{k,l})$, which is
\begin{align*}
 D'_{k,l} = 
\{(x, y, z) \in \ZZ^3 \mid\  & 0 \leq x, y \leq k-1,
(z = -1 \land x + y \leq l-1)\,  \lor \\
&(z = 0 \land l \leq x + y \leq k+l-1 ) \lor\\
&(z = 1 \land k+l \leq x + y) \}.
\end{align*}
Figure 7 depicts digit sets $D'_{k,l}$ for some $k$ and $l$. Note that, 
$D'_{k,l}$ is on the union of three planes $z = -1, z = 0$ and $z = 1 $
if $l > 0$, and 
on the union of two planes $z = 0$ and $z = 1 $
if $l = 0$.
We denote by 
$\SSS_\infty'$,
$\TTT_\infty'$,
$\HHH_\infty'$ the  fractals
$\F(2, D'_{2,0})$, 
$\F(3, D'_{3,0})$, 
$\F(3, D'_{3,1})$, respectively.

Instead of Theorem \ref{mainth1}, we consider the following theorem, which is proved in Section \ref{section:proof}.

\begin{figure}
 \noindent 
\raisebox{1cm}{ $D'_{2,0}$:}
 \includegraphics*[height=28mm]{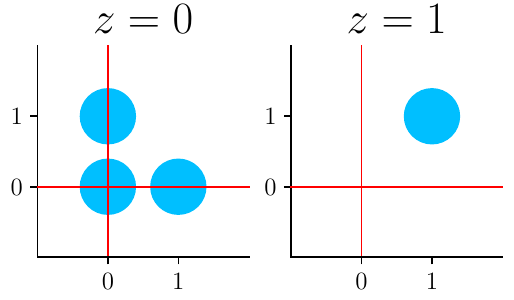} \ \
 \raisebox{1cm}{$D'_{3,1}$:}
  \includegraphics*[height=28mm]{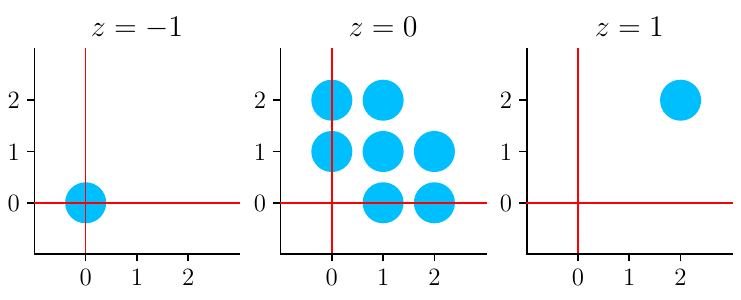}\\
  \raisebox{1cm}{$D'_{3,0}$:}
\includegraphics*[height=28mm]{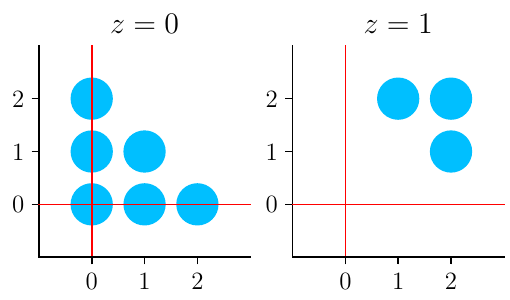} \ \ \ 
\raisebox{1cm}{$D'_{4,1}$:}
\includegraphics*[height=28mm]{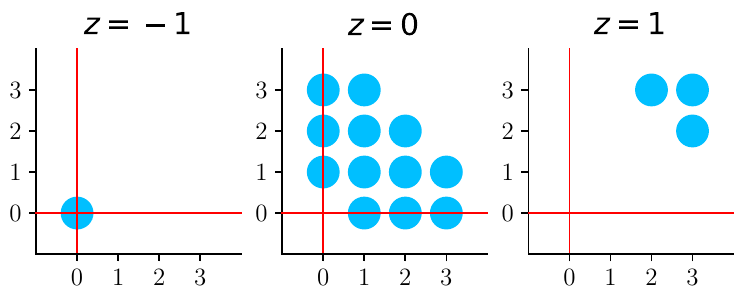}\\
\raisebox{1cm}{$D'_{5,2}$:}
\includegraphics*[height=28mm]{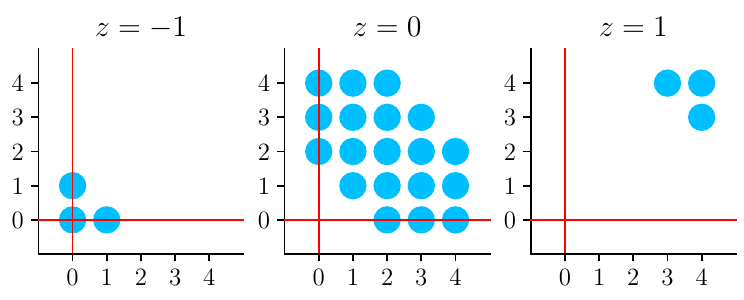}
\caption{\label{fig:digits}
    The digit sets $D_{2,0}'$, $D_{3,0}'$, $D_{3,1}'$, $D'_{4,1}$, and $D'_{5,2}$.
    }
    \end{figure}

    \begin{theorem}\label{mainlemma}
Let $k \geq 2$ and $0 \leq l < \frac{k}{2}$.
$\F(k,D'_{k,l})$ is projected to a set with positive Lebesgue measure if and only if
the projection is done along a coprime vector $(a,  b, c)$ such that
\begin{itemize}
\item[(a)] $a \equiv b \equiv 0\ (\bmod\ k)$  \ (when $(k, l) \ne (3,1)$); 
\item[(b)] $a \equiv b \equiv 0\ (\bmod\ 3)$ or
$a \equiv b \equiv -c\ (\bmod\ 3)$  
\ (when $(k, l) = (3,1)$, i.e., $\F(k, D'_{k,l}) = \HHH'_\infty $).
\end{itemize}
\end{theorem}

We show the equivalence of
Theorem \ref{mainlemma} and Theorem \ref{mainth1}.

\begin{proof}[{\itshape Proof (Equivalence of Theorem \ref{mainlemma} and Theorem \ref{mainth1}).}]
 Because $D'_{k,l}$ is an affine image of $D_{k,l}$ by $\psi_{k,l}$, 
 $\F(k,D_{k,l})$ is projected along $(a,b,c)$  to a set with positive measure if and only if 
 $\F(k,D'_{k,l})$ is projected along 
  $\psi_{k,l}(a,b,c) - \psi_{k,l}(0,0,0) = (a, b, \frac{a+b+c}{k})$ to a set with positive measure.  
 
  Case (a) $(k,l) \ne (3,1)$: We show the following two things. ($\alpha$): 
  If $(a, b, c)$ is a coprime vector such that $k$ and $a + b + c$ are coprime, then for some number $j$, 
  $(p, q, r) = j(a, b,  \frac{a+b+c}{k})$ is a 
 coprime vector such that $p \equiv q \equiv 0\ (\bmod\ k)$.
 ($\beta$): If $(p, q, r)$ is a coprime vector such that $p \equiv q \equiv 0\ (\bmod\ k)$, 
 then for some $m$, $(a, b, c) = m(p, q, kr-p - q)$
 is a coprime vector such that $k$ and $a + b + c$ are coprime.  Statement ($\alpha$) holds for $j = k$ and ($\beta$) holds for $m = \frac{1}{k}$ because $r \not \equiv 0\ (\bmod\ k)$.  
 
 Case (b) $(k,l) = (3,1)$:
 In addition to (a), one can easily see that
   $(a, b, c)$ is a coprime vector such that $a \equiv b  \equiv - \frac{a + b + c}{3} \not \equiv 0\ (\bmod\ 3)$ if and only if
  $(p, q, r) = (a, b,  \frac{a+b+c}{3})$ is a 
 coprime vector such that $p \equiv q \equiv -r\ (\bmod\ 3)$.
\end{proof}

Theorem \ref{mainth1} characterized projections by fixing the arrangement of the fractal and specifying the direction of the projection.
One can also express the characterization
as a property of the projected image as follows.

\begin{corollary}\label{maincor}
Let $d = \frac{1}{k-1}$, $x = ld$, and $y = 1-d$.
On the projected image $X$ of $\F(k, D_{k,l})$, 
consider the four points
$\OO, \PP, \QQ,$ and $\R$ that are the images of
$(x,y,d), (x+d,y,0),(x,y+d,0),and (x+d,y+d,1)$, which are 
elements of $\frac{D_{k,l}}{k-1}$ and thus belong to $\F(k, D_{k,l})$.
The set $X$ has a positive measure if and only if
$r\vec{\OO\R} = p \vec{\OO\PP} + q \vec{\OO\QQ}$ for coprime integers $p, q, r$
such that 
\begin{itemize}
  \item[($a$)]
  $p \equiv q \equiv r \not \equiv 0\ (\bmod\ k)$\  (when $(k, l) \ne (3, 1)$),
  \item[($b$)]
  $p \equiv q \equiv r \not \equiv 0\ (\bmod\ k)$ or
  $p \equiv q \equiv -r \not \equiv 0\ (\bmod\ k)$ 
  \ (when $(k, l) = (3,1)$, i.e., $\F(k, D_{k,l}) = \HHH_\infty $).
  \end{itemize}  
\end{corollary}

\begin{proof}
 Because affine transformations preserve the
 property
 $r\vec{\OO\R} = p \vec{\OO\PP} + q \vec{\OO\QQ}$, 
  we show it for $X$ a projected image of $\F(k,D'_{k,l})$  and 
  $\OO, \PP, \QQ, \R$ the images of
  $(x, y, 0), (x+d, y, 0), (x, y+d, 0), (x+d, y+d, 1)$.
  If the projection is parallel to the $xy$-plane, then $X$ is a null set and 
$r\vec{\OO\R} = p \vec{\OO\PP} + q \vec{\OO\QQ}$ does not hold for $r \ne 0$. Therefore, the statement holds for this case.
 
If not, we can assume that
the projection is along $(a, b, c)$ 
for $c > 0$ to the $xy$-plane and $r > 0$. 
Then,  the points $\OO, \PP, \QQ, \R$ are
$(x, y), (x+d, y), (x, y+d), (x+d-\frac{ad}{c}, y+d - \frac{bd}{c})$,
and $r\vec{\OO\R} = p \vec{\OO\PP} + q \vec{\OO\QQ}$ means $1-\frac{a}{c} = \frac{p}{r}$ and
$1 - \frac{b}{c} = \frac{q}{r}$.
One can see that when $a,b,c$ and $p,q,r$ are coprime integers, 
$a \equiv b \equiv 0\ (\bmod\ k)$ if and only if
$p \equiv q \equiv r \not \equiv 0\ (\bmod\ k)$
because $p = c-a$, $q = c-b$, $r = c$.
One can also see that when $k = 3$, 
$a \equiv b \equiv -c\ (\bmod\ k)$ if and only if
$p \equiv q \equiv -r \not \equiv 0\ (\bmod\ k)$.
\end{proof}

The following result by Kenyon that
appeared as Theorem 14 in \cite{Kenyon1992}  can be derived as a corollary
since an
arbitrary non-colinear four-point set $D \subset \RR^2$ can be affinely transformed to 
\begin{align}
D' = \{(0,0), (1, 0), (0,1), (x, y)\} \label{e:kendigit}
\end{align}
for some $ x, y \in \RR$, that is a projected image of the digit set $D'_{2, 0}$ along the vector $(1-x, 1-y, 1)$.

 \begin{corollary}\label{Cor:mainS}
Let $\D = \{\mathrm{O, P, Q, R}\} \subset \RR^2$ be a non-colinear
four-point set. Then
$\F^2(2, \D)$ has a positive Lebesgue measure if and only if
$r \vec{\mathrm{OR}} = p \vec{\mathrm{OP}} + q \vec{\mathrm{OQ}}$ for odd integers $p, q$, and $r$.
\end{corollary}

\section{Differenced  radix expansion sets}\label{affine}

As stated in the introduction, 
fractals of the form $\F^2(k, D)$ with positive Lebesgue measure
are special cases of self-affine tiles.  A {\em
  self-affine tile} is defined as a compact set $T$ in $\RR^n$ with positive
Lebesgue measure, satisfying $A(T) = \bigcup_{\dd \in D}(T + \dd)$ for an
expansive matrix $A$ and a digit set $D \subset \RR^n$ with cardinality
$\abs{\det(A)}$, and such a tile is usually denoted by $T(A, D)$.
Therefore, $\F^2(k, D)$ 
for $|D| = k^2$ 
with positive measure 
is $T(A, D)$ for $A = \left(\begin{array}{cc}k& 0\\0&k\end{array}\right)$.

Characterization of pairs $(A, D)$ that
generate self-affine tiles
was studied by Bandt \cite{Bandt1991},
Kenyon \cite{Kenyon1992}, Lagarias and Wang \cite{LagariasWang}, and others.
In our study, we use the following 
characterization in \cite{LagariasWang}.
Here, we state it only for self-affine tiles that have the form $\F^n(k,D)$.
\begin{theorem}[Theorem 1.1 of \cite{LagariasWang}]
\label{subtheorem}
  Let $k \geq 2$, $\D \subset \mathbb{R}^n$ be a set of cardinality $k^n$
 and $\zero \in \D$. 
The following four conditions are equivalent.
  \begin{itemize}
  \item[(1)] $\F^n(k,D)$ has positive Lebesgue measure.
  \item[(2)] $\F^n(k,D)$ has nonempty interior.
  \item[(3)] $\F^n(k,D)$ is the closure of its interior, and
its boundary has Lebesgue measure zero.    
  \item[(4)] For each $t \geq 1$, all the 
$(k,D)$-expansions of length $t$
designate distinct points in $\E^n(k, D, t)$, and
$\E^n(k, D)$ is a uniformly discrete set, i.e., 
there exists $\delta > 0$
such that for any distinct elements $\xx, \xx'$ in $\E^n(k, D)$,
$\abs{\xx - \xx'} > \delta$.
\end{itemize}
\end{theorem}

Here, $(k, D)$-expansion, $\E^n(k, D, t) \subset {\mathbb R}^n$, and  $\E^n(k, D) \subset {\mathbb R}^n$ are defined as follows.  
A sequence $(\dd_i)_{0 \leq i < t}$ for $\dd_i \in D$ 
is a {\em $(k, D)$-expansion} of $\xx \in {\mathbb R}^n$ if $\xx$ is expressed as $\xx = \sum_{i=0}^{t-1} k^i\dd_i$.
$\E^n(k, D, t) \subset {\mathbb R}^n$
is the set of points that have $(k, D)$-expansions of length $t$.
It is inductively defined as
\begin{align*}
\E^n(k, D, 1
) &= D\,,\notag\\
\E^n(k, D, t) &= k\E^n(k, D, t-1) + D\,\ (t > 1),
\end{align*}
and we have $\E^n(k, D, 1) \subseteq \E^n(k, D, 2) \subseteq \ldots$.
Finally, the {\em expansion set} $\E^n(k, D) \subset {\mathbb R}^n$ is defined as
\begin{align*}
    \E^n(k, D) &= \bigcup_{t=1}^\infty{\E^n(k, D, t)}\,.%
\end{align*}
In $\E^n(k,D)$, 
we omit $n$ if the dimension is not important. Since $\zero \in D$, 
$\E(k, D)$ satisfies
\begin{align}
    \E(k, D) = k \E(k, D) + D\, . \label{e:eandd}
\end{align}
A digit set all of whose vectors have integer components is called an integral digit set.
\begin{corollary}[Corollary 1.1 of \cite{LagariasWang}] %
    \label{cor-1}
    Let $\D \subset \ZZ^n$ be an integral digit set of cardinality $k^n$ and $\zero \in D$.
    If $\D$ forms a complete residue system of $\ZZ^n/k\ZZ^n$, then
    $\F(k, D)$ has a positive Lebesgue measure.
  \end{corollary}
    \begin{proof}
      We use Theorem \ref{subtheorem}(4) $\to$ (1).
    Suppose that $\D$ forms a complete residue system of $\ZZ^n/k\ZZ^n$.
      $\E^n(k, D)$ is a uniformly discrete set because $D$ is an integral digit set.   If two $(k, D)$-expansions $(\dd_i)_{i < t}$ and 
      $(\ee_i)_{i < t}$ designate the same point, i.e., 
      $\sum_{i=0}^{t-1} k^i \dd_i  = \sum_{i=0}^{t-1}  k^i \ee_i$,
      then $k(\sum_{i=1}^{t-1} k^{i-1} (\dd_i  - \ee_i)) + (\dd_0 - \ee_0) = \zero$ and therefore,
      $\dd_0 - \ee_0 \in k\ZZ^n$ which implies $\dd_0 = \ee_0$.  We can inductively
      show that $\dd_i = \ee_i$ for $0 \leq i < t$.
    \end{proof}
An integral digit set satisfying the condition of Corollary \ref{cor-1} is called a {\em standard digit set}
\cite{LagariasWang1}. It is known that every standard digit set $D$ gives rise to a set $\F(k, D)$ that tiles $\RR^n$ with a lattice tiling \cite{LagariasWang2}.
With this Corollary, one can derive the 'if' part of Theorem \ref{mainlemma}.  That is, 
\begin{lemma}\label{mainlemmapositive}
  Let $a, b, c$ be coprime integers.
  $\F^3(k,D'_{k,l})$ is projected along the vector $(a,  b, c)$  to a set with positive Lebesgue measure if: 
  \begin{itemize}
  \item[(a)] $a \equiv b \equiv 0\ (\bmod\ k)$\  (when $(k, l) \ne (3,1)$); 
   \item[(b)] $a \equiv b \equiv 0\ (\bmod\ 3)$ or
   $a \equiv b \equiv -c\ (\bmod\ 3)$  \ (when $(k, l) = (3,1)$).
  \end{itemize}
\end{lemma}

\begin{proof}
We consider the following two-dimensional digit set 
${D'_{k,l}}^{a,b,c} \subset \ZZ^2$
that is the projection of c$D'_{k,l}$ along $(a, b, c)$.
  \begin{align*}
    {D'_{k,l}}^{a,b,c} = &
    \{(cx + a, cy + b)  \mid\ x, y \in \ZZ ,\  0 \leq x, y \leq k-1,
    \ x + y \leq l-1\}\ \cup \\ 
  &\{(cx, cy) \mid\  x, y \in \ZZ,\ 0 \leq x, y \leq k-1, \ 
  l \leq x + y \leq k+l-1\}\ \cup\\
  &\{(cx - a, cy - b) \mid\ x, y \in \ZZ ,\ 
    0 \leq x, y \leq k-1,
  \ k + l \leq x + y \}.
   \end{align*}
If $a \equiv b \equiv 0\ (\bmod\ k)$, then   
${D'_{k,l}}^{a,b,c}$ is
congruent to 
\[D'' = \{(cx, cy) \mid\  x, y \in \ZZ,\ 0 \leq x, y \leq k-1\}\] modulo
$k\ZZ^2$.
Because $a, b, c$ are coprime and 
$a \equiv b \equiv 0\ (\bmod\ k)$, $c$ is coprime to $k$. Therefore,
$D''$ is a 
complete residue system of $\ZZ^2/k\ZZ^2$.

If $(k,l) = (3,1)$, 
then 
\[{D'_{k,l}}^{a,b,c} = \{(a,b)\} \cup (\{0,1,2\}^2\setminus\{(0,0),(2c,2c)\}) \cup
\{(2c-a,2c-b)\}.\]  It is congruent to $\{0,1,2\}^2$ modulo $3\ZZ^2$ in both cases 
$a \equiv b \equiv -c \equiv 1\ (\bmod\ 3)$ and
$a \equiv b \equiv -c \equiv -1\ (\bmod\ 3)$.
\end{proof}

For an integral digit set $D \subset \ZZ^n$,
a {\em differenced digit set} $\Delta(D)$ is the digit set
\[
\Delta(D) = D + (- D) = \{\xx - \yy \mid \xx, \yy \in D\}
\]
 and
{\em differenced radix expansion set} is the set 
$\E^n(k, \Delta(D))$.
Note that $\zero \in \Delta(D)$ even if $\zero \not \in D$.
Differenced radix expansion sets are important tool for investigating
tiling properties of self-affine tiles \cite{Kenyon1992, LagariasWang}.

To establish the 
`only if' part of Theorem~\ref{mainlemma}, we need to prove the following two assertions:
\begin{enumerate}
\item[1] If coprime vector $(a, b, c)$ fail to meet the conditions specified in the theorem, then the projected image of $\F(k, \D'_{k,l})$ is a null set.
\item[2]If any pair among $a, b$ and $c$ forms an irrational ratio, then the projected image of $\F(k, \D'_{k,l})$ is a null set.
\end{enumerate}
To achieve them, we use Theorem \ref{subtheorem} (1) $\to$ (4), which provides two pathways for establishing that a fractal is a null set: 
by identifying distinct $(k,D)$-expansions for a single point, used to prove assertion 1; and by showing that $\E(k, \varphi(D))$ lacks uniform discreteness, used to prove assertion 2.

For assertion 1,
the existence of distinct $(k,D)$-expansions of length $t$
for a single point means
the presence of   
different sequences $ (\dd_i)_{0 \leq i < t}$ and $ (\ee_i)_{0 \leq i < t}$ from $\D$, satisfying the equation
\[
  \sum_{i=0}^{t-1} k^i \dd_i  = \sum_{i=0}^{t-1}  k^i \ee_i,
\]
or equivalently, 
\[
  \sum_{i=0}^{t-1} k^i (\dd_i  - \ee_i) = \zero.
\] 
This indicates that $\zero$ has a $(k, \Delta(D))$-expansion different from the trivial expansion
$(\zero)_{0 \leq i < t}$.  
Therefore, one can infer that $\F(k, D)$ is a null set by analyzing the differenced radix expansions of the zero vector.

Kenyon's proof of Corollary~\ref{Cor:mainS}, as presented in \cite{Kenyon1992}, employed this  framework. Though it was formulated within the two-dimensional setting, Kenyon's reasoning can be  interpreted as considering the digit set \eqref{e:kendigit} as a projection of a three-dimensional digit set and deriving properties of two-dimensional differenced radix expansion sets through 
the analysis of their three-dimensional counterparts. In the following, we formalize this idea as
a methodology for analyzing projected images of 
fractals of the form $\F^3(k,D)$.

Let $\zero_2$ and $\zero_3$ denote the two and three dimensional zero vectors, respectively. For a three-dimensional digit set $D \subset \RR^3$ and a projection $\varphi$, we have
\[
\Delta(\varphi(\D)) = \varphi(\Delta(\D)).
\]
In addition, 
\[
  \E^2(k, \varphi(\D)) = \varphi(\E^3(k, \D))
\]
holds. Therefore, we have
\[
  \E^2(k, \Delta(\varphi(\D))) = \varphi(\E^3(k, \Delta(\D))).
\]  

Thus, any $(k, \Delta(\varphi(D)))$-expansion of $\zero_2$ is obtained
as a projection of a $(k, \Delta(D))$-expansion of an element in 
$\varphi^{-1}(\zero_2)$.
It means that one can show that $\F^2(k, \varphi(\D))$ is a null set by showing that
$\varphi^{-1}(\zero_2) \cap \E^3(k, \Delta(\D))$ contains an element other than $\zero_3$,
or $\zero_3$ itself has an $(k, \Delta(\D))$-expansion other than the trivial expansion $(\zero_3)_{0 \leq i < t}$.
If $\varphi$ is a projection along $(a, b, c)$, %
then $\varphi^{-1}(\zero_2) = \{(ja, jb, jc) \mid j \in \RR\}$.
Therefore, we have the following theorem.
Note that $\E^2(k, \varphi(D))$ is always a uniformly
discrete set for an integral digit set $D \in \ZZ^3$ and a projection $\varphi$ along
an integral vector.

\begin{theorem}\label{theorem:delta}
    Let $D \subset \ZZ^3$ be an integral digit set of cardinality $k^2$, 
    $a, b, c$ be coprime integers and
    $\varphi$ be a projection along $(a, b, c)$.
    $\F^3(k,D)$ is projected to a null set by $\varphi$ if and only if
    $(ja, jb, jc) \in \E^3(k, \Delta(D))$ for some $j \ne 0$ or
    $\zero_3$ has a $(k, \Delta(D))$-expansion other than 
    those of the form 
    $(\zero_3)_{0 \leq i < t}$.

    \end{theorem}
    
    This theorem allows us to 
    study the directions along which $\F^3(k,D)$ is
    projected to null sets by analyzing the shape of the three-dimensional differenced radix expansion set $\E^3(k, \Delta(D))$.
    Also for assertion 2, 
    analysis of $\E^3(k, \Delta(D))$ proves useful.
    
    \begin{lemma}\label{lemreal}
     $\F^3(k, D)$ is projected 
     to a null set by a projection $\varphi$ if, for all $\delta > 0$, there exists  $\yy \in \E^3(k,\Delta(D))$ such that $0 < \abs{\varphi(\yy)} \leq \delta$.
    \end{lemma}
    \begin{proof}
    We show that  $ \E^2(k, \varphi(D))$ is not uniformly discrete and apply Theorem \ref{subtheorem}.
      For $\E^2(k,\varphi(D))$ not to be  uniformly discrete means
    \begin{align*}
     &\forall \delta > 0.\  \exists {\bm x}, {\bm x}' \in \E^2(k,\varphi(D)). \  \abs{{\bm x} - {\bm x'}} \leq \delta \land \xx \ne \xx'\\
     \Leftrightarrow\ &\forall \delta > 0.\  \exists {\bm x} \in \E^2(k,\Delta(\varphi(D))). \  0 <  \abs{\bm x} \leq \delta \\ %
     \Leftrightarrow\ &\forall \delta > 0.\  \exists {\bm x} \in \varphi(\E^3(k,\Delta(D))). \  0 < \abs{\bm x} \leq \delta \\ %
     \Leftrightarrow\ &\forall \delta > 0.\  \exists {\bm y} \in \E^3(k,\Delta(D)). \ 0 < \abs{\varphi({\bm y})} \leq \delta %
    \end{align*}
    \end{proof}
    
    Though Theorem \ref{theorem:delta} and Lemma \ref{lemreal} are applicable to any integral digit set $D$, the structure of $\E^3(k, \Delta(D))$  appears to be complicated in general, making it challenging to derive useful results directly from them. This is why we focus our analysis on layered fractal imaginary cubes $\F^3(k, D_{k,l})$ and introduce affine-transformed fractal objects $\F^3(k, D'_{k,l})$. In the next section, we will analyze the structure of $\E^3(k, \Delta(D'_{k,l}))$ to derive the `only if' parts of Theorem \ref{mainlemma}.

    \section{Proofs of the main theorem}
    \label{section:proof}
    
    We prove the `only if' part of Theorem~\ref{mainlemma}.
    We first consider the case of $S_\infty = \F(2, D'_{2,0})$, for which we have a simplified proof. Subsequently, we present a proof applicable to $\F(k, D'_{k,l})$ in general.
    
\subsection{Sierpinski tetrahedron case}\label{sec:sier}

\begin{figure}
\centering
 \includegraphics*[height=30mm]{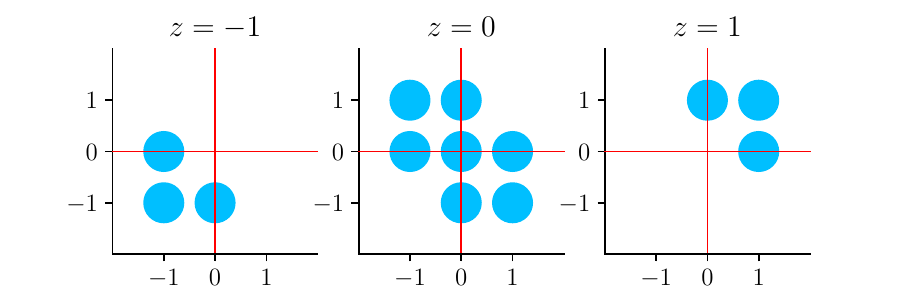}
\caption{\label{fig:dds}
The differenced digit set $\Delta(D'_{2,0})$. }
\end{figure}

\begin{figure}  
\includegraphics*[width=140mm]{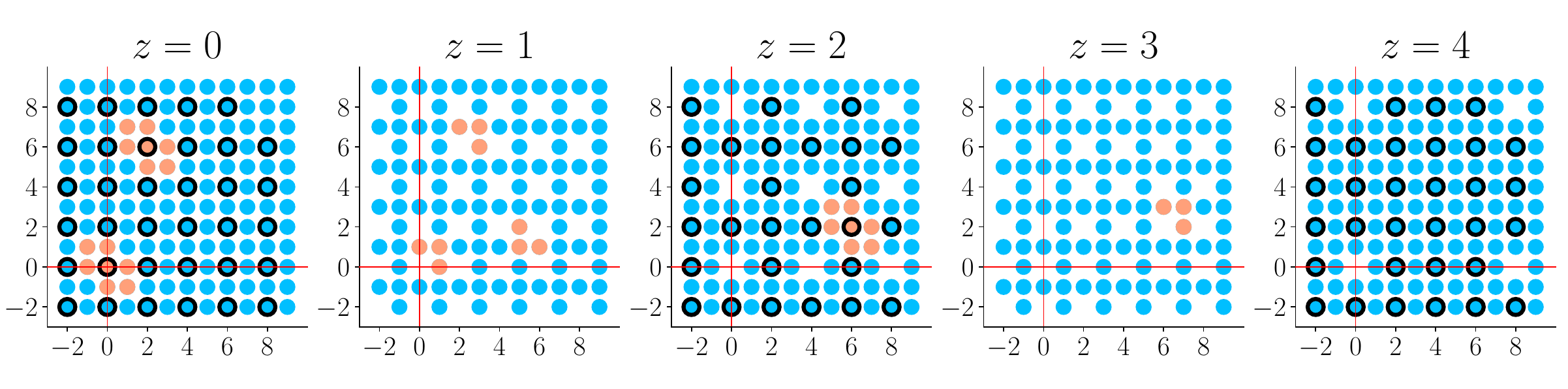}  
\caption{\label{fig:des}
Points in $\E(2, \Delta(\D_{2,0}'))$ for $-2 \leq x,y \leq 9$,
$0 \leq z \leq 4$  are represented by blue and orange filled circles. Among them, points 
in $2\E(2, \Delta(D_{2,0}'))$  have black borders and
translations of $\Delta(D_{2,0}')$ by $(0,0,0)$, $(2,6,0)$ and $(6,2,2)$ are
colored orange.
}
\end{figure}

We analyze the projected images of $\SSS'_\infty$ by examining the structure of $\E(2,\Delta(D'_{2,0}))$. For a subset $A \subseteq \ZZ^3$ and an integer $z$, we use the notation $\s{A}{z}$ to represent the slice $\{(x, y) \mid (x, y, z) \in A\}$. 
We simplify notation by introducing $\DD$ as a shorthand for $\Delta(D'_{2,0})$ and $\EE$ as a shorthand for $\E(2,\Delta(\D'_{2,0}))$.
The differenced digit set $\DD$ is decomposed into slices as follows:
\begin{align*}
    \DDs{1} &= \{(x,y) \mid x, y \leq 1 \land x + y \geq 1\},\\
  \DDs{0}\ \  &= \{(x,y) \mid \abs{x}, \abs{y}, \abs{x + y} \leq 1\},\\
  \DDs{-1} &= \{(x,y) \mid x, y \geq -1 \land  x + y  \leq -1\}.
\end{align*}
These sets are depicted in Figure~\ref{fig:dds}. 
By equation \eqref{e:eandd},
the slices of the differenced radix expansion set $\EE$ 
satisfies the followings:
\begin{align}
 \EEs{2n} &= 2 \EEs{n} + \DDs{0}\, ,\label{l:s2}\\
 \EEs{2n+1} &= 2 \EEs{n+1} + \DDs{-1}\ \cup\ 2 \EEs{n} + \DDs{1}\, \label{l:s3}.
\end{align}
Some parts of $\EE$ are calculated and depicted in Figure \ref{fig:des}.
Note that it is not an inductive definition because 
$\E_{0}$ depends on $\E_{0}$, and $\E_{1}$ depends on $\E_{0}$ and $\E_{1}$. 
Note also that \eqref{l:s3} has the form of the union of two sets, which makes the 
analysis of the structure of $\EE$ complicated.
Contrary to this,  we will see in Lemma \ref{lem:ds} below that $\EE$ has a rather simple form. To show this, 
we define the set %
$C \subset \ZZ^3$ as follows.
\begin{align*}
B &=  2\ZZ^2\, ,\\
\s{C}{0} &=  \ZZ^2\, ,\\
\s{C}{2^mc} &=  \ZZ^2 \setminus 2^mB\ \ (m\geq 0, c \text{ is odd}).
\end{align*} 

\begin{lemma}\label{lem:ds}
  $ \EE = C$.
 \end{lemma}
  \begin{proof}
     
 $ \EE \subseteq   C$ 
  is demonstrated by proving 
that for all $t$, 
$\E(2, \DD, t) \subseteq C$  holds by  induction on $t$, that is straightforward and we omit it.
 $C \subseteq  \EE $ is demonstrated by proving  
\[
 \forall  {\xx} \in \ZZ^3,   {\xx} \in C \Rightarrow   {\xx} \in  \EE 
 \]
 through well-founded induction on 
  $| {\xx}|$, where $\abs{(x,y,z)} =  \abs{x} + \abs{y} + \abs{z}$.
   It holds for the case $\abs{ {\xx}} = 0$
  because $\zero_3 \in C$ and $\zero_3 \in  \EE$.
 For $ {\xx} = (x, y, z) \in C$ satisfying $\abs{ {\xx}} > 0$, 
 we define ${\xd } \in C$ such that
 $\xx \in 2\xd + \DD$ as follows.

  Case $z = 0$ (Refer to $z=0$ column of Figure 9). If $(x,y) \in \DDs{0}$, we set 
  ${ \xd} = \zero_3$.
If $(x,y) \not \in \DDs{0}$, 
the family of sets $\{2(x', y') + \DDs{0}  \mid (x', y') \in \ZZ^2\}$ covers  $\ZZ^2$.  Therefore,
 we can choose ${\xd} = (x', y', 0) \in \ZZ^2 \times \{0\} \subset C$
 such that $(x, y) \in 
 2(x', y') + \DDs{0}$. 

 Case $z = 1$ (Refer to $z=1$ column of Figure 9).
 The family $\{2(x', y') + \DDs{1} \mid (x', y') \in \ZZ^2\}$ covers  $\ZZ^2 \setminus B = \s{C}{1}$ without overlaps.
 Therefore, we 
 set ${\xd}$ to be
 the unique $(x', y', 0)$ such that $(x, y) \in 2(x', y') + \DDs{1}$.

 Case $z = -1$: Similar.

Case $z$ is even and $z \ne 0$ (Refer to $z=2$ column of Figure 9).
 If $y$ and $z$ are even, then we set  ${\xd} = (x/2, y/2, z/2) \in C$.
Otherwise, %
 the family $\{2(x', y') + \DDs{0} \mid (x', y') \in \ZZ^2\}$ double covers $\ZZ^2 \setminus B$.
 On the other hand, the family 
 $\{2(x', y') + \DDs{0} \mid (x', y') \in B\}$ is disjoint.  Therefore, 
 \[\{2(x', y') + \DDs{0} \mid (x', y') \in \ZZ^2 \setminus B\}\]
 covers $\ZZ^2 \setminus B$.
 Thus, we can choose ${\xd} = (x', y', z/2) \in C$
 such that 
\[(x, y) \in 
 2(x', y') + \DDs{0}.\] 

 Case $z$ is odd and $|z| \ne 1$ (Refer to $z=3$ column of Figure 9). 
 In this case, $(x, y) \not \in B$.
 The family 
\[\{2(x', y') + \DDs{1} \mid (x', y') \in \ZZ^2\} \cup
 \{2(x', y') + \DDs{-1} \mid (x', y') \in \ZZ^2\}\]
 is a double cover for $\ZZ^2 \setminus B$.
 On the other hand, the family 
 \[\{2(x', y') + \DDs{1} \mid (x', y') \in B\} \cup
 \{2(x', y') + \DDs{-1} \mid (x', y') \in B\}\]
 is disjoint.
 Therefore, the family
 \[\{2(x', y') + \DDs{1} \mid (x', y') \in \ZZ^2 \setminus B\} \cup
 \{2(x', y') + \DDs{-1} \mid (x', y') \in \ZZ^2 \setminus B\}\] covers $\ZZ^2 \setminus B$.
 Since $\s{C}{z'} \supseteq \ZZ^2 \setminus B$ for 
 $z' = \frac{z-1}{2}$ and $z' = \frac{z+1}{2}$, we can choose
 $\xd  = (x', y', \frac{z-1}{2}) \in C$ such that
 $(x,y) \in 2(x',y') + \DDs{1}$ 
 or
 $\xd  = (x', y', \frac{z+1}{2}) \in C$ such that
 $(x,y) \in 2(x',y') + \DDs{-1}$.
  
 One can see that in all the cases, $\abs{{\xd}} < \abs{{\xx}}$.
 Therefore, $\xd \in \EE$ by induction hypothesis. 
 Thus, $\xx \in 2\xd + \DD \subset 2 \E + \DD = \E$.
  \end{proof}

  {\em Note:} As this proof illustrates, formulating a $(2, \DD)$-expansion for a given vector $(x, y, z) \in C$ would become complex. This complexity necessitates a case analysis based on congruence modulo 4 for both $x$ and $y$ coordinates when $z$ is odd.

\medbreak

Now, we prove the 'only if' part of Theorem \ref{mainlemma}
for the case $k = 2, l = 0$.
\begin{lemma}\label{lemmaS1}
If $(a, b, c)$ is a coprime vector such that either $a$ or $b$ is an odd number, then 
$\SSS'_\infty$ is projected along $(a,b,c)$ to a null set.
    
\end{lemma}  

\begin{proof}
By Lemma \ref{lem:ds}, we know that 
$(a, b, c) \in \E$.
Thus, the statement follows from Theorem \ref{theorem:delta} with $j = 1$. 
\end{proof}

\begin{lemma}\label{lemma:S2} 
    If any pair among $a, b$ and $c$ forms an irrational ratio, then 
    $\SSS'_\infty$ is projected along the vector $(a,b,c)$ to a null set.
 \end{lemma}
\begin{proof}
  In the case where $c = 0$, let $\varphi$ be the projection along $(a,b,c)$.  Since $\varphi$ maps the three points $(0,0,0)$, $(1, 0, 0)$, and $(0, 1, 0) \in D'_{2,0}$  
  onto the same line, $ \E(2, \varphi(D'_{2,0}))$ is not uniformly discrete. The statement follows from Theorem \ref{subtheorem}.
  
Suppose that $c \ne 0$. 
Let $u = \frac{a}{c}$ and  $v = \frac{b}{c}$. 
We consider the projection $\varphi$ along $(a,b,c)$ to the $xy$-plane, whose image is
$\F(2, D_\SSS^{u,v})$ for
\[D_\SSS^{u,v} = \{(0,0), (1, 0), (0, 1), (1-u, 1-v)\}.\] 
 
Let $\lfloor x \rfloor$ and $\{x\}$ be the
integer and fractional parts of $x \in \RR$, respectively.  For arbitrary $j$, 
$\abs{\varphi(\lfloor ju \rfloor, \lfloor jv \rfloor, j)} =
\abs{\varphi(ju - \{ju\}, jv - \{jv\}, j)} = \abs{(\{ju\},
  \{jv\})}$ because  
$\varphi(u, v, 1) = 0$.
We show the following claim and
apply Lemma \ref{lemreal}.
\begin{quote}
{\bf Claim.}
For all $\delta > 0$, there  exists $j \in 2\ZZ$ such that \[(\lfloor ju \rfloor, \lfloor jv \rfloor, j) \in \E 
\land 0 < |(\{ju\}, \{jv\})| < \delta.\]
\end{quote}
  
We first consider the case when $u$ is irrational and $v = \frac{p}{q}$ is rational. The set
$\{\{j'u\}  \mid j' \in q\ZZ\}$ is a dense subset of $[0,1]$.
Therefore, there exists some $j'$ such that
 $0 < \{j'u\}- \frac{1}{2} < \frac{\delta}{2}$ and
 $\{j'v\} = 0$.  Since
  $0 < 2\{j'u\}- 1 < \delta$, we obtain that 
$0 < 2(j'u-\lfloor j'u \rfloor) - 1 < \delta$. 
Thus, for $j = 2j'$, 
$0 < ju-2\lfloor j'u \rfloor - 1 < \delta$.
  It means that $\lfloor ju \rfloor = 2\lfloor j'u \rfloor + 1$  and 
 $0 < \{ju\} < \delta$, and therefore
  $(\lfloor ju \rfloor, \lfloor jv \rfloor, j) \in C$ and
 $0 < |(\{ju\}, \{jv\})| < \delta$.
 If $u$ and $v$ are both irrational,  
 we have a similar argument starting with a dense subset
 $\{(\{j'u\}, \{j'v\})  \mid j' \in \ZZ\}$ of $[0,1]^2$ and  
 choosing $j'$ such that  $0 < \{j'u\}- \frac{1}{2} < \frac{\delta}{4}$ and
  $0 < \{j'v\}- \frac{1}{2} < \frac{\delta}{4}$.
  \end{proof}

\subsection{The general case}

\begin{figure}[t]
    \begin{flushleft}   
      \raisebox{1cm}{$\Delta(D'_{3,0}):$ } \hspace*{0.9cm}
      \includegraphics*[width=75mm]{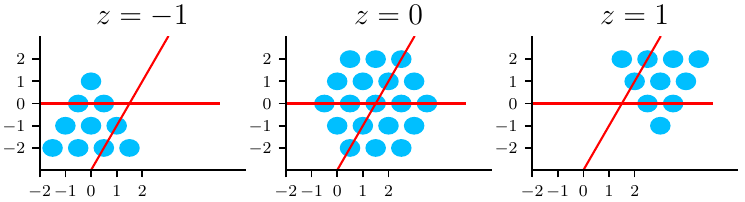} \\
      \raisebox{1cm}{$\Delta(D'_{3,1}):$}
     \includegraphics*[width=125mm]{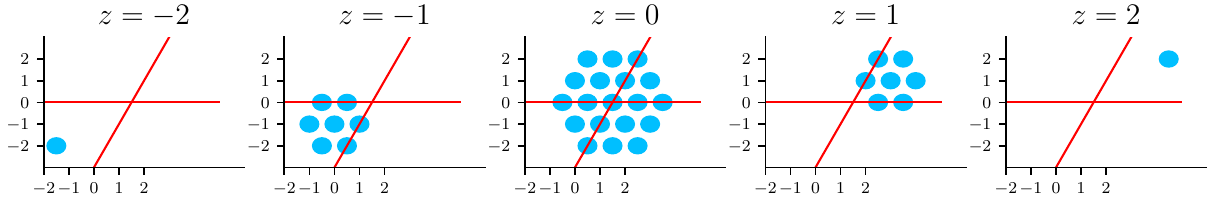}\\
     \raisebox{1cm}{$\Delta(D'_{4,1}):$}
     \includegraphics*[width=125mm]{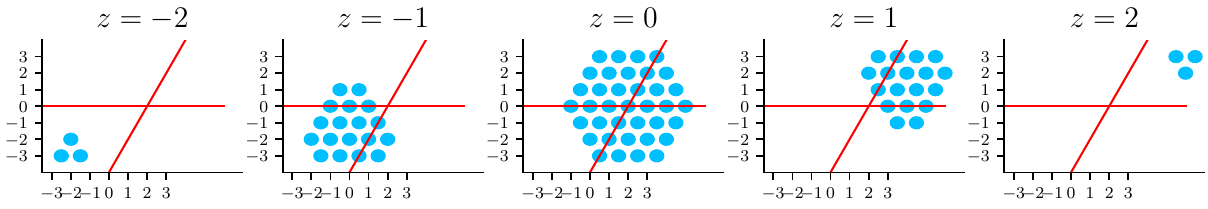}\\
     \raisebox{1cm}{$\Delta(D'_{5,2}):$}
     \includegraphics*[width=125mm]{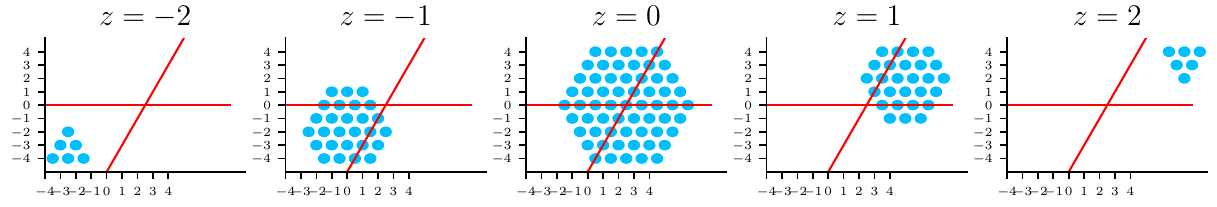}
    \end{flushleft}
     \caption{\label{fig:ddigits}
    The differenced digit sets $\Delta(D'_{3,0})$, $\Delta(D'_{3,1})$,  $\Delta(D'_{4,1})$ and $\Delta(D'_{5,2})$. In this and the following pictures, 
    the elementary vectors $(1, 0)$ and $(\frac{1}{2}, \frac{\sqrt{3}}{2})$ are used to
    demonstrate their three-fold symmetries.}
    \end{figure}

We proceed with the proof of the `only if' part of Theorem \ref{mainlemma} for all cases $2 \leq k$ and $0 \leq l < \frac{k}{2}$, including the specific instances of $\TTT'_\infty = \F(3, \D'_{3,0})$ and $\HHH'_\infty = \F(3, \D'_{3,1})$. 
In this subsection, we fix $k$ and $l$ and 
denote $\Delta(D'_{k,l})$ and $\E(k
,\Delta(\D'_{k,l}))$ as $\DD$ and $\EE$, respectively. Note that, though we simplify our notation by omitting parameters $k$ and $l$, 
these sets, along with others introduced subsequently such as $B$, $B'$, $G(m)$, $G'(m)$ and $A(m)$, depend on the chosen values of $k$ and $l$. 
Slices of $\DD$ have the following forms. 
Here, $\DDs{-2}$ and $\DDs{2}$ are empty sets when $l = 0$.  
\begin{align*}
  \DDs{2} &= 
  \{(x, y) \in \ZZ^3 \mid\  |x|, |y| \leq k-1,\  k + 1 \leq x + y \}\ \ \  (\text{when }l > 0)\\
  \DDs{1} &= 
  \{(x, y) \in \ZZ^3 \mid\  l-k+2 \leq x, y\leq k-1,\  1 \leq x + y \leq 2k-l-2\}  \\
  \DDs{0} &= 
  \{(x, y) \in \ZZ^3 \mid\  |x|, |y| \leq k-1,\  \abs{x + y} \leq k-1\}  \\
  \DDs{-1} &= 
  \{(x, y) \in \ZZ^3 \mid\  1-k \leq x, y \leq k-l-2,\  -2k+l+2\leq x + y \leq -1 \}  \\
 \DDs{-2} &= 
 \{(x, y) \in \ZZ^3 \mid\  |x|, |y| \leq  k-1,\   x + y \leq -k - 1   \}\ \ \   (\text{when }l > 0)
\end{align*}
These sets are depicted in Figure~\ref{fig:ddigits} for the cases $(k, l) = (3, 0)$,
$(3, 1)$, $(4, 1)$ and $(5,2)$.
In these and subsequent figures, the elementary vectors 
$(1,0)$ and $(\frac{1}{2}, \frac{\sqrt{3}}{2})$
are used to demonstrate their three-fold symmetries after the corresponding linear transformation.
We sometimes use
the hexagonal norm on $\RR^2$ defined as 
\begin{align}
H(x,y) = \max(|x|, |y|, |x+y|)\,. \label{e:hex}
\end{align}
For every $k$ and $l$,  $\Delta_0$ is the set of points 
in the disk of radius $k-1$ with this norm, that is a regular hexagon with these elementary vectors.
Refer to $z = 0$ column of Figure \ref{fig:ddigits}.

If $(k, l)$ is not $(3, 1)$ nor $(4, 1)$,
the following relation holds 
by equation \eqref{e:eandd}. 
\begin{align*}
    \E_{kn+i} &= k \E_{n} + \Delta_{i}\hspace*{0.5cm} (i = -2, -1, 0, 1, 2)\, 
\end{align*}
When $(k, l) = (4, 1)$,
$\Delta_{2}$ and $\Delta_{-2}$ are not empty and $kn+2 = k(n+1)-2$ holds.
Therefore,
in this case, 
\begin{align*}
  \E_{4n+i} &= 4 \E_{n} + \Delta_{i}\hspace*{0.5cm} (i =  -1, 0, 1),\\
  \E_{4n+2} &= 4 \E_{n} + \Delta_{2} \cup 4 \E_{n+1} + \Delta_{-2}\, .
 \end{align*}
Overlaps also occur when  $(k, l) = (3, 1)$, and we have
\begin{align}
  \E_{3n-1} &= 3 \E_{n} + \Delta_{-1} \cup 
  3 \E_{n-1} + \Delta_{2}\, ,\label{line:eh4}\\
  \E_{3n} &= 3 \E_{n} + \Delta_{0}\, ,\nonumber\\
  \E_{3n+1} &= 3 \E_{n} + \Delta_{1} \cup 
  3 \E_{n+1} + \Delta_{-2}\, . \nonumber
 \end{align}

\begin{figure}[t]
 \begin{flushleft} 
  $\E(3,\Delta(D'_{3, 0})):$    \\
  \includegraphics*[width=140mm]{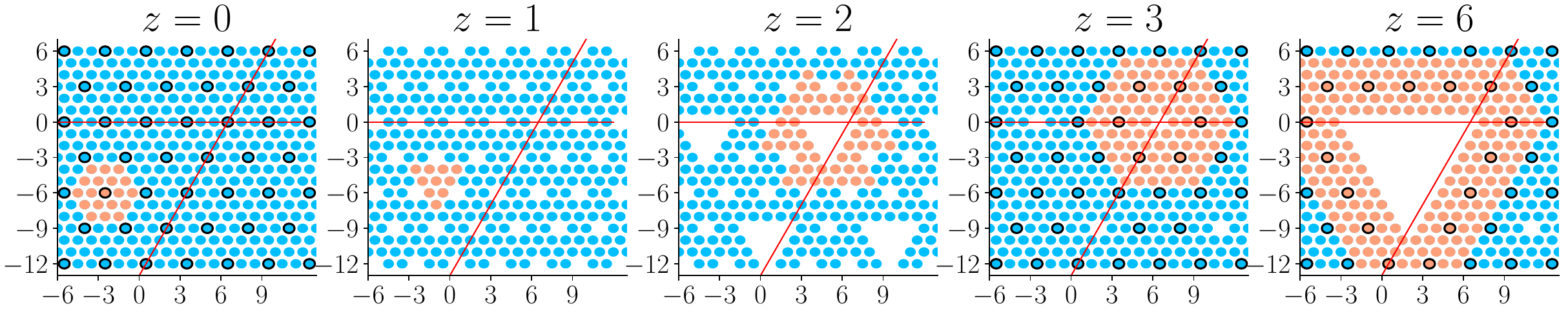}\\
  $\E(3,\Delta(D'_{3, 1})):$    \\
  \includegraphics*[width=140mm]{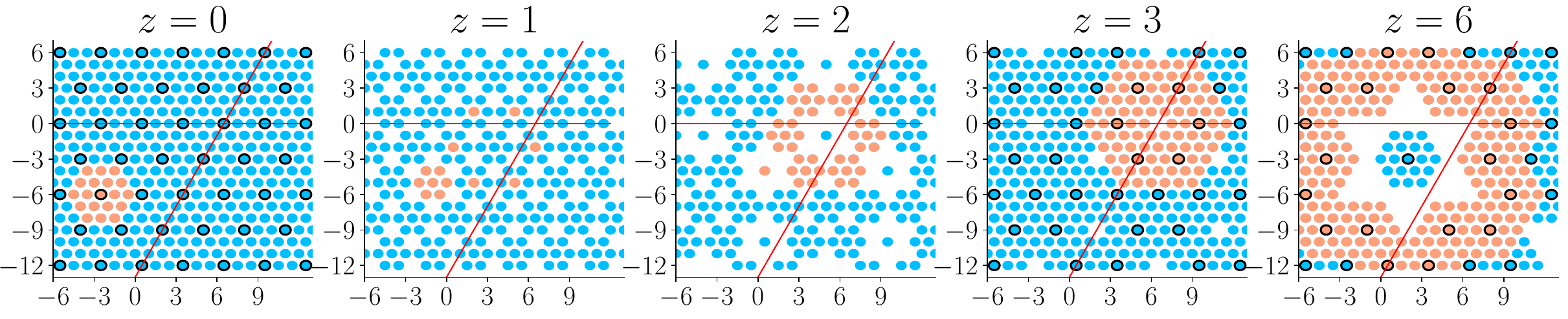}\\
  $\E(5,\Delta(D'_{5, 2})):$    \\
  \includegraphics*[width=140mm]{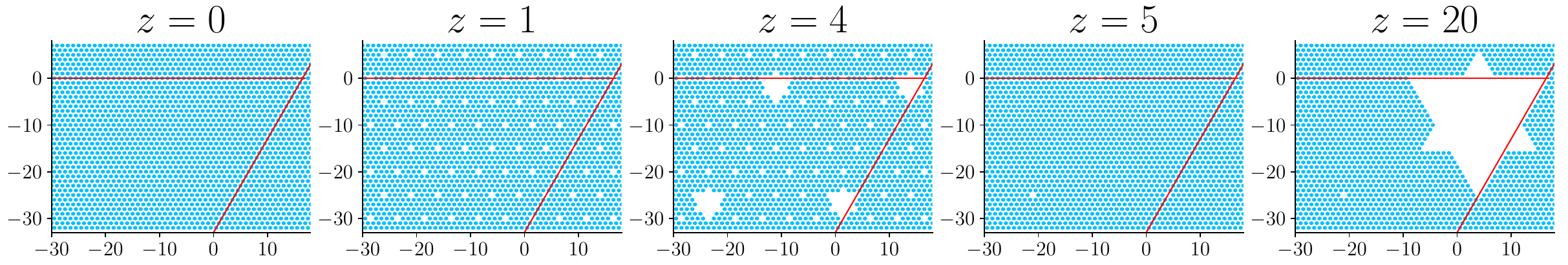}\\ 
 \end{flushleft}
  \caption{\label{fig:pt5} 
 Slices at $0, 1, k-1, k, k(k-1)$ of
 $\E(3, \Delta(\D_{3,0}'))$,
 $\E(3, \Delta(\D_{3,1}'))$, and
 $\E(5, \Delta(\D_{5,2}'))$.
 In the first two lines of figures,
  points in $3\E$ have black borders and 
  $\xx + \Delta$ for some $\xx \in 3\E$ are colored orange.
 }
 \end{figure}

The structures of $\E$ for $k > 2$ are quite complicated
and one could not easily characterize them as we did for the case $(k, l) = (2, 0)$.  However, 
for the following proof, we only use the slices
$\E_{k^m(k^n-1)}$ for $m, n \geq 0$, that have rather simple and uniform structure as we will see in
Lemma \ref{lem-E}.
$\E_z$ for $z = 0, 1, k-1, k, (k-1)k$ are depicted in Figure \ref{fig:pt5} for the cases $(k, l)$ is $(3, 0)$, $(3, 1)$, and $(5, 2)$.

In Lemma~\ref{lem-E}, we express
$\E_{k^m(k^n-1)}$ using the sets $B$, $B'$, and $A(m)$ ($m = 0, 1, \ldots$) defined below.
The sets $B$ and $B'$ are defined as follows:
\begin{align*}
B &= \left\{\begin{array}{ll}
 k\ZZ^2 & (\text{$(k,l) \ne (3,1)$})\\
 k\ZZ^2  \cup  k\ZZ^2 + (1,1)  & (\text{$(k, l) = (3, 1)$}) 
\end{array}\right.\\
B' &= \left\{\begin{array}{ll}
k\ZZ^2 & (\text{$(k,l) \ne (3,1)$})\\
k\ZZ^2  \cup  k\ZZ^2 + (-1,-1)  & (\text{$(k, l) = (3, 1)$}) 
\end{array}\right.
\end{align*}    
As we are going to see, $B$ is the complement of $\E_1$ and $B'$ is the complement of $\E_{-1}$. Thus, for a visual representation of $B$, refer to the $z=1$ column in Figure \ref{fig:pt5}.
We have
\begin{lemma}\label{lem:16}
\begin{itemize}
\item[(1)] $k\ZZ^2 + \DD_{1} = \ZZ^2 \setminus B'$.  
\item[(2)] $k\ZZ^2 + \DD_{-1} = \ZZ^2 \setminus B$. 
\end{itemize}
\end{lemma}
\begin{proof}
    (1) When $(k,l) \ne (3,1)$, the set
    $\DD_1$ is congruent to $\{0,1,\ldots,k-1\}^2\setminus \{(0,0)\}$ modulo $k\ZZ^2$.
   When $(k,l) = (3,1)$,
   $\DD_1$ is equal to $\{0,1,2\}^2\setminus \{(0,0), (2,2) \}$.
   The proof of (2) is similar. 
   \end{proof}

The sets 
$A(m) (m = 0,1,\ldots)$ are defined as follows:
\begin{align*}
G(m) &= \{(x,y) \in \ZZ^2\mid x,y \leq 0 \land x + y \geq
-k^{m+1} \}\\
G'(m) &= \{(x, y) \in \ZZ^2\mid x, y \geq -\frac{(k+1)k^m }{2}\land x + y \leq
 -\frac{(k-1)k^m}{2} \}\\
A(m) &= \left\{\begin{array}{ll}  
G(m) & (\text{except for the followings})\\
G(m) \cup G'(m)  
&   (\text{$k$ is odd and $l = \frac{k-1}{2}$}) \, \\
\end{array}\right.
 \end{align*}
 $A(m)$ represents the void around the origin in $\E_{k^m(k^n-1)}$ for any $n > 0$. 
 In the third and fifth columns of Figure \ref{fig:pt5}, $A(0)$ and $A(1)$ appear as void spaces around the origin, for the cases of $(k, l)$ being $(3,0)$, $(3,1)$, and $(5,2)$. Please note that the central 'island' within 
 these voids, observable in the case of $(k,l) = (3,1)$, should be disregarded.
 
Key properties of $A(m)$ are given in Lemma 
\ref{lemforlem-E}.
For that purpose, we
 introduce some notations:
 \begin{align*}
    V & = \{(0,0), (-k,0), (0,-k)\}\, \\
    E & = \{(-1, 1), (0, 1), (1, 0), (1,-1), (0, -1), (-1,0)\}\, \\
    \Nei(X) & = (X + E) \setminus X\, 
\end{align*}
$E$ is on the unit circle with respect to the hexagonal norm  \eqref{e:hex}, and
$\Nei(X)$ denotes the set of points in the outer boundary layer of $X$, i.e., each located one unit away from the nearest point in X with respect to the hexagonal norm.

 \begin{lemma}\label{lemforlem-E}
  \begin{enumerate}
 \item[(1)]$\{(0, 0)\} = \DD_0 \setminus 
(kE + \DD_0)$. 
\item[(2)]  
$A(0) = (\DD_{-1} \setminus  (kE + \DD_{-1})) \cup V$.  
\item[(3)]  $A(m+1) = (kA(m) + \DD_{0}) \setminus
(kN(A(m)) + \DD_{0})$. 
\end{enumerate}
\end{lemma}
\begin{proof}

\noindent
(1) Easily checked.

\noindent
(2) We prove $A(0) = (\Delta_{-1} \cup V) \setminus (kE + \Delta_{-1})$.
The six copies, $\{k\xx + \Delta_{-1} \mid \xx \in E\}$, 
of $\Delta_{-1}$ collectively create a void space around $(0,0)$. 
We show that it coincides with $A(0)$.
Considering that $\Delta_{-1}$ possesses a hexagonal or triangular outline, the void typically manifests as a 'star-like' or a triangular configuration. For $l \leq \frac{k}{2} - 1$, the edges of $k(-1, 1) + \DD_{-1}$ and $k(0, 1) + \DD_{-1}$ merge, forming a continuous edge. This pattern of connection repeats symmetrically for the pairs ($k(1, 0) + \DD_{-1}$, $k(1, -1) + \DD_{-1}$) and ($k(0, -1) + \DD_{-1}$, $k(-1,0) + \DD_{-1}$), yielding a triangular void represented by $G(0)$.
Refer to $z = 2$ slice of $\E(3,\Delta(D'_{3,0}))$ in Figure~\ref{fig:pt5} for an example.
The other case $l > \frac{k}{2} - 1$ is exactly the case $k$ is odd and $l = \frac{k-1}{2}$ because
$\frac{k}{2}  > l$.  
Under these circumstances, each of the aforementioned edge pairs exhibits a single-point separation, culminating in a 'star-like' formation identified as $G(0) \cup G'(0)$.
Refer to $z = 2$ and $4$ slices of $\E(3,\Delta(D'_{3,1}))$ and
$\E(5,\Delta(D'_{5,2}))$,
respectively,
in Figure~\ref{fig:pt5} for examples.

(3) We examine the structure of the void space around $(0,0)$
formed by the set $\{k\xx + \DD_{0} \mid \xx \in \Nei(A(m))\}$, where each copy of $\DD_{0}$ is placed $k$ units away from one another and from the outer edge of the set $kA(m)$ with respect to the hexagonal norm.  
It is important to note that
$\DD_{0}$ is larger than $\DD_{-1}$, resulting in no gaps between edges of adjacent elements as observed in (2). 
Consequently, this arrangement preserves the overall shape of 
$kA(m)$, effectively creating what we identify as 
$A(m+1)$. Refer to the $z = 6$ slices of $\E(3,\Delta(D'_{3,0}))$ and
$\E(3,\Delta(D'_{3,1}))$, and the $z = 20$ slice of $\E(5,\Delta(D'_{5,2}))$
in Figure~\ref{fig:pt5} for examples.
\end{proof}

With these preparations, we show the necessary properties of some of the slices of ${\mathcal E}$.

\begin{lemma}\label{lem-E}
\begin{enumerate}
\item[(1)]   $\E_{0} = \ZZ^2$.
\item[(2)]   $\E_{k^n} = \ZZ^2 \setminus k^{n}B'$ for $n \geq 0$.
\item[(3)]  $\E_{k^m(k^n-1)} \supseteq \ZZ^2 \setminus k^{m}B \setminus  (k^{m+n}B'+A(m))$  for $m\geq 0$, $n > 0$.
\end{enumerate}
\begin{proof}
\noindent
(1)   
Adopting the same approach as the proof for the $z = 0$ case presented in Lemma \ref{lem:ds}.

\noindent
(2)  Induction on $n$. 

Case $n = 0$: (Refer to $z = 1$ column of Figure~\ref{fig:pt5}.)
We show that
$\E_1 = \ZZ^2 \setminus B'$.
When $(k, l) \ne (3, 1)$,
$\E_1 = k\E_0 + \Delta_1 =
k \ZZ^2 + \Delta_1$,
which is equal to
$\ZZ^2 \setminus B'$ by 
Lemma \ref{lem:16}(1).
When $(k, l) = (3, 1)$,
$\E_1 = 3\E_0 + \Delta_{1} \cup 3\E_1 + \Delta_{-2}$ and
$3\E_1 + \Delta_{-2} \subseteq 3\ZZ^2 + (-2, -2) \subseteq 3\ZZ^2 + \Delta_1$.  Therefore, also for this case, 
$\E_1$ is equal to 
$\ZZ^2 \setminus B'$.%

Case $n > 0$:(Refer to $z = k$ column of Figure~\ref{fig:pt5}.)
We assume that
$\E_{k^{n-1}} = \ZZ^2 \setminus k^{n-1}B'$ and 
show that $\E_{k^{n}} = \ZZ^2 \setminus k^{n}B'$.  We have
\begin{align*}
\E_{k^{n}} &= k\E_{k^{n-1}} + \DD_0 = 
k(\ZZ^2 \setminus k^{n-1}B') + \DD_0\\
&= (k\ZZ^2 \setminus k^{n}B') + \DD_0\,.
\end{align*}
Since $k\ZZ^2 + \DD_0 = \ZZ^2$  
and no points of $kE$ belong to $k^nB'$,
it follows that this is equal to $\ZZ^2 \setminus (k^nB' + K)$
where $K$ is the subset of $\DD_0$ that does not intersect with  $\xx + \DD_0$ for every $\xx \in kE$.
That is, $K = \DD_0 \setminus  (kE + \DD_0)$, which is $\{(0, 0)\}$
by Lemma \ref{lemforlem-E}(1).
Thus, $\EE_{k^{n}} = \ZZ^2 \setminus k^{n}B'$.

\noindent
(3) Induction on $m$.  

Case $m = 0$ 
(Refer to $z = k-1$ column of Figure~\ref{fig:pt5}).
We prove the inclusion 
\[\EE_{k^n-1} \supseteq \ZZ^2 \setminus B \setminus (k^nB' + A(0))\] for $n > 0$.  
In both $(k, l) = (3, 1)$ and $(k, l) \ne (3, 1)$ cases, we have  that
$\EE_{k^{n}-1} \supseteq k\EE_{k^{n-1}} + \Delta_{-1}$. Therefore, 
\begin{align*}
  \EE_{k^{n}-1} 
    &\supseteq k\EE_{k^{n-1}} + \Delta_{-1}
  = k(\ZZ^2 \setminus k^{n-1}B') + \DD_{-1} \\
  &= (k\ZZ^2 \setminus k^{n}B') + \DD_{-1}.
\end{align*}
We have $k\ZZ^2 + \DD_{-1} = \ZZ^2 \setminus B$ by Lemma \ref{lem:16}(2).  In addition,  no points of $kE$ belong to $k^nB'$.
Therefore, 
it is equal to $(\ZZ^2\setminus B) \setminus (k^nB' + K)$
for $K = \Delta_{-1} \setminus (kE + \Delta_{-1})$ as in (2).
By Lemma \ref{lemforlem-E}(2),
$K$ becomes
$A(0)$ by adding the three points $(0,0), (-k,0), (0,-k)$, which are already contained in $B$. 
  Thus, \[\EE_{k^{n}-1} \supseteq \ZZ^2 \setminus B \setminus (k^nB' + A(0)).\]

  Case $m > 0$: By induction hypothesis
  \begin{align*}
  \EE_{k^m(k^n-1)} &= 
  k \EE_{k^{m-1}(k^n-1)}  + \DD_{0}\\
  &\supseteq k(\ZZ^2 \setminus k^{m-1}B \setminus (k^{m-1+n}B'+A({m-1}))) + \DD_{0}\\
  &= k\ZZ^2  \setminus k^{m}B \setminus (k^{m+n}B'+kA(m-1)) + \DD_{0}\,.
  \end{align*}
  Since $k\ZZ^2 + \DD_0 = \ZZ^2$,  it is equal to $\ZZ^2 \setminus (k^mB +K)\setminus (k^{m+n}B' + K')$ where
  $K$ is $\{(0,0)\}$ as in (2), and
  $K'$ is the subset of $kA(m-1) + \DD_0$
that do not intersect with  $\xx + \DD_0$ for $\xx \in k\ZZ^2
\setminus kA(m-1)$.
Note that $k\Nei(A(m-1))$ coincides with the set of points 
$\xx \in k\mathbb{Z}^2 \setminus kA(m-1)$ such that $\xx + \DD_0$ intersects with $kA(m-1) + \DD_0$, because $\DD_0$ is the disk of radius $k-1$ with hexagonal norm.
Therefore, $K'$ is the subset of $kA(m-1) + \DD_0$
that do not intersect with  $\xx + \DD_0$ for $\xx \in 
k\Nei(A(m-1))$, that is, 
  \begin{align*}
   K' & = (kA(m-1) + \DD_{0}) \setminus (k\Nei(A(m-1)) + \DD_0)
  \end{align*}
  which is $A(m)$
   by Lemma \ref{lemforlem-E}(3).
   Thus, \[\EE_{k^{m}(k^n-1)} \supseteq \ZZ^2 \setminus k^{m}B \setminus (k^{m+n}B' + A(m)).\]
\end{proof}
\end{lemma}
Note: In Lemma \ref{lem-E}(3), the specific shape of $\E_{3^m(3^n-1)}$ was not detailed as it was not necessary for our discussion. The equality in this lemma holds when $(k,l) \ne (3, 1)$. 
When $(k,l) = (3, 1)$, 
a precise specification of $\E_{3^m(3^n-1)}$ can be given: within $\E_{3^m(3^n-1)}$, each $A(m)$-shaped void contains a hexagonal 'island'.
This island takes the form of the intersection of $G(m)$ and $G'(m)$, which 
is then reduced by removing points along its boundary. This fact
is proved inductively: First, 
\[
\E_{3^n-1} = 3 \E_{3^{n-1}} + \Delta_{-1} \cup 
  3 \E_{3^{n-1}-1} + \Delta_{2} 
\]  
by \eqref{line:eh4}.
This formula allows for the inductive specification of $\E_{3^n-1}$ based on the knowledge of $\E_{3^{n-1}}$. Consequently, $\E_{3^m(3^n-1)}$ is specified by induction on $m$ according to the formula
\[
\E_{3^m(3^n-1)} = 3\E_{3^{m-1}(3^n-1)} + \Delta_0.
\] 

\medbreak

Now, we prove the 'only if' 
part of Theorem \ref{mainlemma}.
In the proof of Lemma \ref{lemmaS1},  we showed that
$(a,b,c) \in \E^3(2, 
\Delta(D_{2,0}'))$
and then applied Theorem \ref{theorem:delta} for $j = 1$.
However, in the proof of Lemma \ref{lemmakl1} below,  
$(a,b,c) \not \in \E^3(k, \Delta(D_{k,l}'))$ in general and therefore the same argument does not apply. 
Instead, we prove that 
\[(ja, jb, jc) \in \E^3(k, \Delta(D_{k,l}'))\] for some $j > 0$.

\begin{lemma}\label{lemmakl1}
Let $k \geq 2$ and $0 \leq l < \frac{k}{2}$. 
Suppose that $a, b, c$ are coprime integers such that 
\begin{enumerate}
\item[(a)] $\neg(a \equiv b \equiv 0 \ (\bmod\ k))$ \ (when $(k,l) \ne (3,1)$);
\item[(b)] $\neg(a \equiv b \equiv 0\ (\bmod\ 3)) \land
\neg (a \equiv b \equiv -c \ (\bmod\ 3))$\ (when $(k,l) = (3,1)$).
\end{enumerate}
Then, $\F(k, D'_{k,l})$ is projected along $(a, b, c)$ to a null set.
\end{lemma}

\begin{proof}
First, we construct numbers $s, m, n$ and a sequence $t_0, t_1,\ldots$ such that 
\[
t_p s c = k^m(k^{2^pn}-1)
\]
 for every $p  \geq 0$ as follows.
Consider an $s'$ for which $s' c$ has the form $k^m c'$, where $c'$ is an integer that is coprime to $k$.
We choose the smallest such $s'$, that is not a multiple of $k$ nor shares any common factors with $c'$.
By applying Euler's theorem, $k^{\phi{(c')}} \equiv 1\pmod{c'}$
for $\phi$ the Euler's totient function.
It follows that, for $n = \phi{(c')}$ and some $s'' > 0$, 
$s''c' = k^n - 1$. Consequently, for $s = s's''$, we have
$sc = k^m(k^n - 1)$.  
Let $t_{0} = 1$ and 
$t_p = (k^{2^{p-1}n} + 1) \ldots (k^{2n} + 1)(k^n + 1)$ for $p > 0$.  
They satisfy  $t_1 s c = k^m (k^n + 1)(k^n - 1) = k^m (k^{2n} - 1)$. Similarly, we can show by induction that $t_p s c = k^m (k^{2^p n} - 1)$ for all $p \geq 0$.

Next, we show the following claim:
\begin{quote}
{\bf Claim~1.} $(t_psa, t_psb) \not \in  k^mB$.
\end{quote}
Case $m > 0$:
Suppose for contradiction that $(t_psa, t_psb) \in  k^mB$.
Since $t_p$ and $s''$ are coprime to $k$, 
both $s'a$ and $s'b$ are multiples of $k^m$.
Since $s'c$ is also a multiple of $k^m$ and
$s'$ is not a multiple of $k^m$,
it follows that $a, b, c$ share a common factor, which contradicts the assumption that $a, b, c$ are coprime.

Case $m = 0$: In this case, $s' = 1$.
Therefore, the condition $(a, b) \not \in k\ZZ^2$ implies 
$(t_psa, t_psb) \not \in  k\ZZ^2$.
This proves the case $(k,l) \ne (3,1)$.
When $(k,l) = (3,1)$, 
we also need to show $(t_psa, t_psb) \not \in  3\ZZ^2 + (1,1)$.
In this case, either $a$ or $b$ is not congruent to $-c$ modulo 3.
Therefore, $t_psa$ or $t_psb$ is not congruent to $-t_psc$ modulo 3.
Since $t_psc = 3^n-1 \equiv -1\ (\bmod\ 3)$, it follows that  $t_psa$ or $t_psb$ is not congruent to 1 modulo 3.

Now, we fix $p$ such that $2^{p} n \geq m+1$.
We show the following claim:
\begin{quote}
  {\textbf Claim~2.} $(ja, jb, jc) \in \E$ holds for at least one of $j = t_ps$ or $j = t_{p+1}s$.
  \end{quote}
Let $h = k^{m+{2^p} n+1}$ and $C = (k^{{2^p} n} + 1)A({m})$.
For the sake of contradiction, suppose that both 
$(t_psa, t_psb, t_psc) \not \in \E$ and
$(t_{p+1}sa, t_{p+1}sb,t_{p+1}sc) \not \in \E$ hold.
Since $t_{p'}sc = k^m(k^{2^{p'}}-1)$, using Lemma 19(3) and
Claim~1, we can show that for $p' \geq 0$, the condition $(t_{p'}sa, t_{p'}sb, t_{p'}sc) \not \in \E$ implies $(t_{p'}sa, t_{p'}sb) \in  k^{m+2^{p'}n} B' 
+A({m})$.
Therefore, we have the following.
\begin{align*}
(t_psa, t_psb) &\in  k^{m+2^{p}n}B'+A({m})\,\\
(t_{p+1}sa, t_{p +1}sb) &\in  k^{m+2^{p+1}n}B'+A(m)\,
\end{align*}
From the first formula, we have
\begin{align*}
(t_{p+1}sa, t_{p +1}sb) &\in 
(k^{{2^p} n} + 1)(k^{m+2^{p}n}B'+A({m})) 
\end{align*}
We show that the two sets
\[(k^{{2^p} n} + 1)(k^{m+2^{p}n}B'+A({m}))
\text{ and }k^{m+2^{p+1}n}B'+A(m)\] 
do not intersect outside of $k^mB$,
which is a contradiction since $(t_{p+1}sa, t_{p +1}sb) \not \in k^mB$ by Claim 1.

Case (a) $(k,l) \ne (3,1)$:
We have
\begin{align}
    (k^{{2^p} n} + 1)(k^{m+2^{p}n}B'+A({m})) 
= & (k^{{2^p} n} + 1)(k^{m+2^{p}n+1}\ZZ^2 +A({m})) \nonumber\\
= & (k^{{2^p} n} + 1)h\ZZ^2+C   \label{e:1st}
\end{align}
and 
\begin{align}
k^{m+2^{p+1}n}B'+A(m) 
=  k^{2^{p}n}h\ZZ^2+A(m)\,.   \label{e:2nd}
\end{align}

Modulo $h\ZZ^2$, the sets 
\eqref{e:1st} and \eqref{e:2nd} are congruent to $A(m)$
and $C$, respectively.
For $q = (k^{2^pn} + 1)k^m$,
$C$ is contained in the triangle with  vertices at
$\{(0,0), (-kq, 0), \linebreak (0, -kq)\}$ 
when $k$ is even or $l \ne \frac{k-1}{2}$,
and $C$ is contained in the union of this triangle and another
one with vertices at
\[\{(-\frac{(k+1)q}{2}, q), (q,-\frac{(k+1)q}{2}), (-\frac{(k+1)q}{2}, -\frac{(k+1)q}{2})\}\] 
when $k$ is odd and $l = \frac{k-1}{2}$.
Refer to  Figure \ref{fig:figCA}(a) which depicts the points in $C$ and some copies of $A(m)$ modulo $h\ZZ^2$ for
the case $(k,l) = (5,2)$ with the parameters $m= 0,n=1,p=0$.
Clearly, the second triangle does not intersect with copies of $A(m)$ modulo $h\ZZ^2$. 
Since $kq = h + k^{m+1}$, the first triangle intersects with copies of $A(m)$ modulo $h\ZZ^2$ at
the three vertices.  With the hexagonal norm, 
$k^{2^pn} + 1$, which is the distance between any two adjacent points of $C$, is greater than $k^{m+1}$, the distance between 
pairs of points among the vertices 
$\{(0,0), (0,-k^{m+1}), (-k^{m+1},0)\}$
of $A(m)$,
because $2^pn \geq m+1$. 
Therefore, these three vertices are the only points shared by $C$ and copies of $A(m)$ modulo $h\ZZ^2$. 
However, these three points belong to $k^mB$. 
Therefore, the two sets defined in \eqref{e:1st} and \eqref{e:2nd}
do not intersect outside of $k^mB$.

\begin{figure}[t]
    \centering
  \subfloat[$k=5$, $l=2$, $m=p=0$, $n=1$. Thus $h=25$ and $C=6A(0)$ on the slice $t_{p+1}sc=24$.]
  {\includegraphics*[height=40mm]{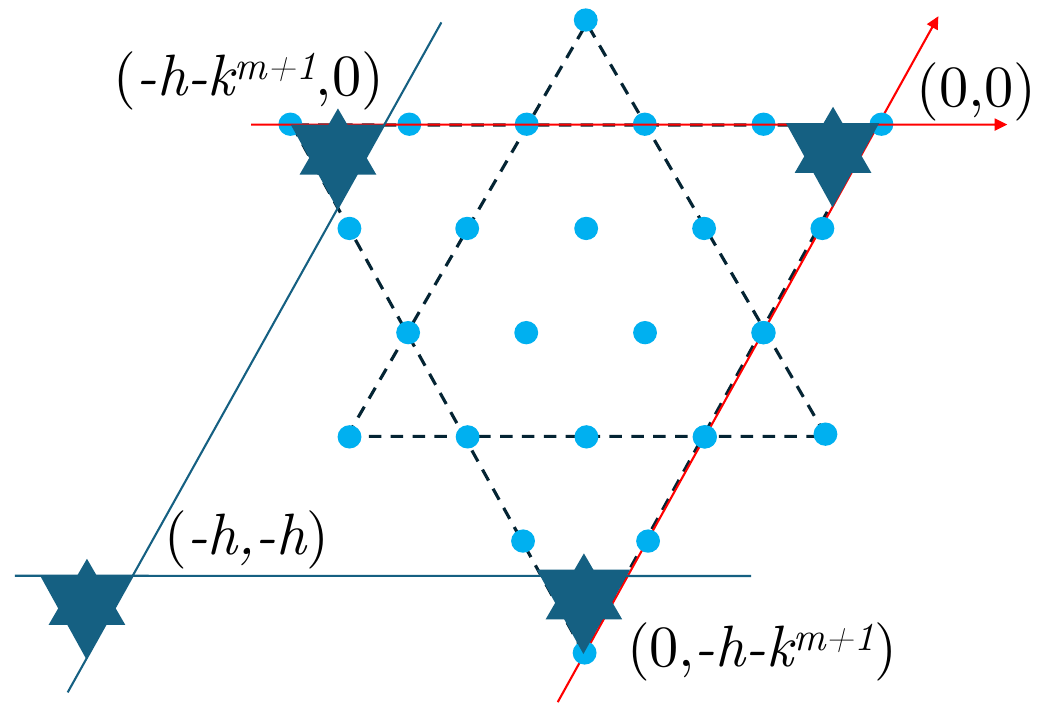}}\qquad
    \subfloat[$k=3$, $l=1$, $m=p=0$, $n=1$. Thus $h=9$ and $C=4A(0)$ on the slice $t_{p+1}sc=8$.]{\includegraphics*[height=40mm]{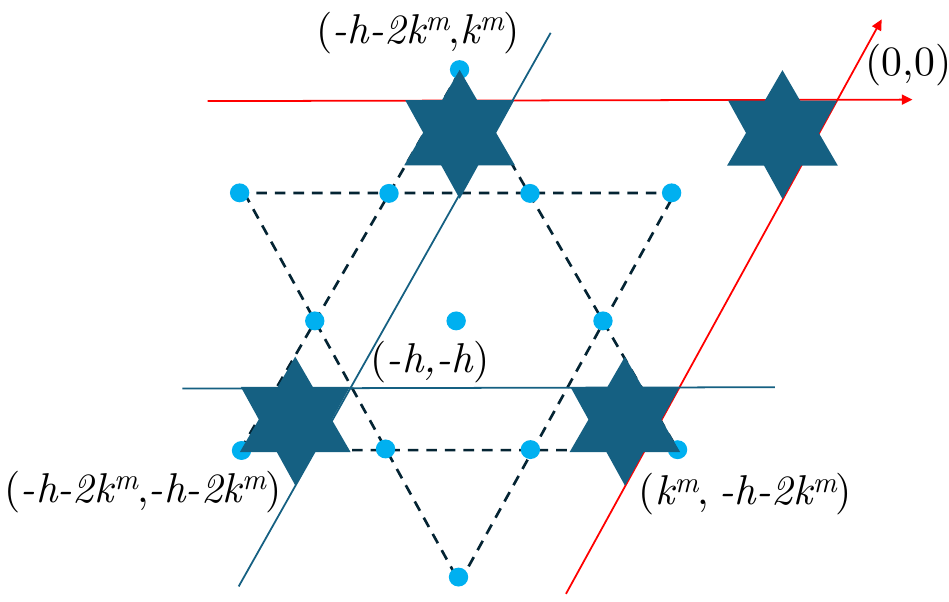}} 
   \caption{\label{fig:figCA}  
      Configurations of copies of A(m) and C.      
          Points in $C$ are shown as blue points, and copies of $A(0)$ modulo $h\ZZ$ are located within dark 'starlike' regions. In (b), only the second copy of $C$ is depicted. 
    }
   \end{figure}

Case (b) $(k,l) = (3,1)$: 
Instead of \eqref{e:1st} and \eqref{e:2nd}, we have the following:
\begin{align}
&(k^{{2^p} n} + 1)(k^{m+2^{p}n}B'
+A({m})) \nonumber \\
= &(k^{{2^p} n} + 1)((k^{m+2^{p}n+1}\ZZ^2 +A({m})) \cup
(k^{m+2^{p}n}(k\ZZ^2+(-1,-1)) +A({m}))) \nonumber\\
= & ((k^{{2^p} n} + 1)h\ZZ^2+C)  \cup 
((k^{{2^p} n} + 1)h\ZZ^2  + (k^{{2^p} n} + 1)k^{m+2^{p}n}(-1,-1) +C  ) \label{e:1st2}
\end{align}
and
\begin{align}
&k^{m+2^{p+1}n}B'
+A(m)  \nonumber\\ 
= & (k^{2^pn}h\ZZ^2+A(m)) 
\cup  (k^{m+2^{p+1}n}(k\ZZ^2+(-1,-1))+A(m)).  \label{e:2nd2}
\end{align}
The set \eqref{e:2nd2} is congruent to $A(m)$ modulo $h\ZZ^2$, and the set \eqref{e:1st2} is congruent to
$C \cup (-\frac{h}{3},-\frac{h}{3})+ C$ modulo $h\ZZ^2$, because
$k^{m+2^{p}n} = \frac{h}{3}$ for the case $k = 3$.
$C$ does not intersect with copies of $A(m)$ modulo $h\ZZ^2$ outside of $k^mB$, as we studied in (a).
$(-\frac{h}{3},-\frac{h}{3})+ C$
consists of two triangles that are translations of those
we studied for the case (a) by $(-\frac{h}{3},-\frac{h}{3})$.
Refer to  Figure \ref{fig:figCA}%
(b) which depicts 
$(-\frac{h}{3},-\frac{h}{3})+ C$ and some copies of $A(m)$ modulo $h\ZZ^2$ for the case $(k,l) = (3,1)$ with the parameters $m= 0,n=1,p=0$.
Clearly, the first triangle does not intersect with copies of $A(m)$. 
For the second one, the three vertices are translations by $(-\frac{h}{3},-\frac{h}{3})$ of $\{(-2q, q), (q,-2q), (-2q, -2q)\}$ 
for $q = k^{m+2^pn} +k^m = \frac{h}{3} + k^m$,
which  are \[\{(-h - 2k^m, k^m), (k^m, -h - 2k^m), (-h - 2k^m, -h - 2k^m)\}.\]
These three points are congruent to $(- 2k^m, k^m)$, $(k^m, - 2k^m)$ and $(- 2k^m, - 2k^m)$, respectively,  modulo $h\ZZ^2$, 
and they also belong to $A(m)$ modulo $h\ZZ^2$.
These three vertices are the only points shared by 
$(-\frac{h}{3},-\frac{h}{3})+ C$ and copies of $A(m)$ modulo 
$h\ZZ^2$ as we studied for the case (a).
However, they belong to $k^mB$.
Therefore, the two sets defined in \eqref{e:1st2} and \eqref{e:2nd2}
do not intersect outside of $k^mB$.
This finishes the proof of Claim~2.

By Claim~2, 
$(ja, jb, jc) \in \E$ for some $j > 0$. Therefore, by Theorem~\ref{theorem:delta}, 
$\F(k, D'_{k,l})$ is projected along $(a, b, c)$ to a null set.
\end{proof}

Finally, we prove the irrational part of Theorem \ref{mainlemma}.

\begin{lemma}\label{lemmaT2}
  If any pair among $a, b$ and $c$ forms an irrational ratio, then 
  $\F(k,D'_{k,l})$ is projected along the vector $(a,b,c)$ to a null set.
  \end{lemma}
  \begin{proof}
  If $c = 0$, let $\varphi$ be the projection along $(a,b,c)$.
  Since the digits on the plane $z = 0$, 
  whose cardinality is at least $k+1$,
  are mapped by $\varphi$ onto the same line, $ \E(k, \varphi(D'_{k,l}))$ is not uniformly discrete. The statement follows from Theorem \ref{subtheorem}.
  
  If $c \ne 0$, let $u = \frac{a}{c}$, $v = \frac{b}{c}$ and $\uu = (u,v)$. 
  We consider the projection $\varphi$ to the $xy$-plane and aim to prove 
   \begin{align*}
  \forall r>0. \exists \yy \in \E. \  0 < \abs{\varphi(\bm y)} \leq \frac{1}{k^r}\,
  \end{align*}
  in order to apply Lemma \ref{lemreal}.   In Lemma \ref{lemma:S2}, the point  $\yy$ was chosen from even slices of $\E(2, \Delta'_{2, 0})$.
  On the other hand, we need to choose $\yy$ from the $k^m(k^n-1)$-th slices of  $\E(k, \Delta'_{k, l})$ for
  $m \geq 0, n > 0$ as we %
  understand the structure of $\E$ only on these  slices.
  Moreover, since these slices contain large voids, it is crucial  to avoid them when selecting $\yy$. 
  Consequently, our proof methodology here differs from that used in Lemma \ref{lemma:S2}.

  For $\bm x = (x, y) \in \RR^2$, we define $\lfloor \bm x \rfloor$ as $(\lfloor x \rfloor, \lfloor y \rfloor)$ and $\{ \bm x \}$ as $(\{ x\}, \{ y\})$.
  For a given $r > 0$, we partition the unit square $[0,1)^2$ into $k^{2r}$ disjoint regions and set
  $f(s, t) = [\frac{s}{k^r}, \frac{s+1}{k^r})\times 
  [\frac{t}{k^r}, \frac{t+1}{k^r})$ for $0 \leq s, t < k^r$.
  By the pigeonhole principle, at least one of these regions, say  $f(s, t)$, contains an infinite number of points from  
  $\{\{k^{i} \bm u \} \mid i \geq 0\}$.
  Let $0 \leq g_0 <  g_1 < \ldots$ be the enumeration of 
  $\{i \mid \{k^{i} \bm u \} \in f(s,t)\}$.
  We then define the numbers $h$, $m$, and $n_i$ for $i = 0, 1, \ldots$ as follows: let 
  $h = g_0$, $m = g_1$, and $n_i = g_{i+2} - m$.

  Let $\bm x_i = \lfloor k^{m+n_i} \bm u \rfloor - \lfloor k^{m} \bm u \rfloor$,
  $(x_i, y_i) = \bm x_{i}$, and
  $\yy_i = (x_{i}, y_{i}, k^{m}(k^{n_i} - 1))$.
  For each $i$, $\yy_{i}$ 
    satisfies $\abs{\varphi(\yy_{i}
    )} \leq \frac{1}{k^r}$ with respect to the maximum norm because
  \begin{align*}
  \abs{\varphi(\yy_{i}
  )} &= \abs{(x_{i},  y_{i}) - k^{m}(k^{n_i} - 1) \bm u}\\
  &   = \abs{\xx_i - k^{m+n_i} \bm u + k^{m}\bm u} \\
  & = \abs{ - \{ k^{m+n_i} \bm u \} + \{ k^{m}\bm u\}} < \frac{1}{k^r}.
  \end{align*}
  
  Therefore, we complete the proof by showing that 
  $\yy_{i} \in \E$ for some $i$.
  By Lemma \ref{lem-E}, it means  
    $\bm x_i \in \ZZ^2 \setminus k^{m}B \setminus
    (k^{m+n_i}B'+A(m))$.
  We demonstrate the following claim
  \begin{quote}
    {\bf Claim.}  For some $i$,  both $\bm x_i \not \in k^{m}B$ and 
  $\bm x_i \not \in k^{m+n_i}B' + A(m)$ hold.
  \end{quote}
  Suppose, for the sake of contradiction, that 
  $\bm x_i \in k^{m}B$ or $\bm x_i \in k^{m+n_i}B' + A(m)$ for every $i$.  If
  $\bm x_i = \lfloor k^{m+n_i} \bm u \rfloor - \lfloor k^{m} \bm u \rfloor \in k^{m+n_i}B' + A(m)$
  then
  \[
    k^{m+n_i} \bm u - k^{m} \bm u \in k^{m+n_i}B' + A(m) + [-\frac{1}{k^r},
  \frac{1}{k^r}]^2.
  \]
  Therefore, 
  \begin{align}
    \bm u &\in \frac{k^{n_i}}{k^{n_i} -1}B'  + \frac{1}{k^{n_i}-1}
  \frac{A(m)+[-\frac{1}{k^r},\frac{1}{k^r}]^2}{k^{m}}\notag\\
    &\subseteq \frac{k^{n_i}}{k^{n_i} -1}B'  + \frac{1}{k^{n_i}-1}[-k-1,1]^2.\label{line-41}
  \end{align}
  If this holds for arbitrary large $n_i$, then we have $\bm u \in B'$.
  However, it
  contradicts  the fact that at least one
  component of $\bm u$ is irrational.
  Therefore, $\bm x_i \in k^{m+n_i}B' + A(m)$ only for a finite number of $i$. 
  Let $I$ be their maximum.
  Thus,  $\bm x_i \in k^{m}B$ for $i > I$.
  
  Let $(u_j)_{j \geq 0}$ and $(v_j)_{j \geq 0}$ be the $k$-ary expansions of 
  the fractional parts of $u$ and $v$, respectively,  and
  let $\sigma, \tau \in \{0,1,\ldots,k-1\}^r$ be the $k$-ary expressions of  $s$ and $t$, respectively. We say that an infinite sequence $(u_j)_{j \geq 0}$ contains a sequence $\sigma$ 
  at position $n$ if, for $r$ the length of $\sigma$,  $(u_j)_{n \leq j < n+r}$ is equal to $\sigma$.  According to the definition
  of the indices $h, m, n_i\ (i = 0,1,\ldots)$, 
  $(u_j)_{j \geq 0}$ and $(v_j)_{j \geq 0}$ contain $\sigma$ and 
  $\tau$, respectively, at $h, m$ and $m+n_i$ ($i = 0,1,\ldots$).
  
  Since $\bm x_i = \lfloor k^{m+n_i} \bm u \rfloor - \lfloor k^{m} \bm u \rfloor
   \in k^{m}B \subset k^m\ZZ^2$ for $i > I$, the last $m$ digits of the integral parts of $k$-ary expansions of $k^{m}u$
  and $k^{m+n_i}u$  coincide for $i > I$.  
  Let $\alpha$ be this sequence of length $m$.
  Then, $(u_j)_{j \geq 0}$ contains $\alpha$ at the indices 0 and $n_i$ $(i > I)$.
  Since $(u_j)_{j \geq 0}$ contains $\sigma$ at %
  $m$ and $n_i+m$ $(i > I)$,
  $(u_j)_{j \geq 0}$ contains the sequence $\alpha\sigma$ at 0 and $n_i$ for $i > I$,
  On the other hand, the sequence $\alpha\sigma$ contains $\sigma$ twice at indices $h$ and $m$.
  Therefore, $(u_j)_{j \geq 0}$ contains $\sigma$ at $n_i+h$ and $n_i+m$.
  The same is true for $(v_j)_{j \geq 0}$ and $\tau$, and $(v_j)_{j \geq 0}$ contains $\tau$ at $n_i+h$ and $n_i+m$.
  Since both $n_i+h$ and $n_{i-1}+m$ are the last index before $n_{i}+m$
  at which $\sigma$ and $\tau$ are contained in $(u_j)_{0 \leq j}$ and $(v_j)_{0 \leq j}$ simultaneously,
  $n_i + h = n_{i-1} +m$ holds. Thus, we have $n_i - n_{i-1} = m -  h$ and the subsequence $(u_j)_{n_{i-1} \leq j < n_i}$
  is equal to the sequence $\alpha_{0 \leq j < m-h}$ for all $i > I$.  Thus, $u$ is rational, %
and similarly $v$ is  also rational. %
   This contradicts  the fact that $u$ or $v$ is irrational. 
  This completes the proof of the Claim.
  \end{proof}
   
\section{Concluding remarks}\label{sec:conclusion}

We have characterized the directions along which layered fractal imaginary cubes, including the Sierpinski tetrahedron,
T-fractal, and H-fractal, are projected to sets with positive measures. As explained in Theorem \ref{mainth1}, such directions are uniformly described, with the exception of the H-fractal case.
The H-fractal possesses a unique property that
it exhibits six-fold symmetry, which effectively 
'doubles' its projected images with positive measures.

A natural question would be whether our result can be generalized to fractal imaginary cubes,
that is, imaginary cubes of the form $\F^3(k, D)$ with $|D| = k^2$, in general.  About fractal imaginary cubes,  the following proposition holds immediately.

\begin{figure}
\centering
\subfloat[A  digit set $D$ expressed as a Latin square. Here, $\mathrm {Lat}(L) = \{(x, y, L(x,y)) \mid x, y \in \ZZ, 0 \leq x, y \leq 3 \}$.]{
\begin{minipage}[b]{5cm}
\[
D = \mathrm{Lat}\left(\begin{bmatrix}
  1& 3 & 0 & 2\\
  2&0 & 1 & 3\\
  3& 1& 2 & 0\\
  0& 2& 3 & 1
\end{bmatrix}\right)\\
\qquad
\]

\end{minipage} }\qquad
\subfloat[The projected image of
$\F(4, D)$ along $(1,1,1)$. It is positively measured by Corollary \ref{cor-1}.]
{\qquad \raisebox{-0.9cm}{\includegraphics[width=3.5cm]{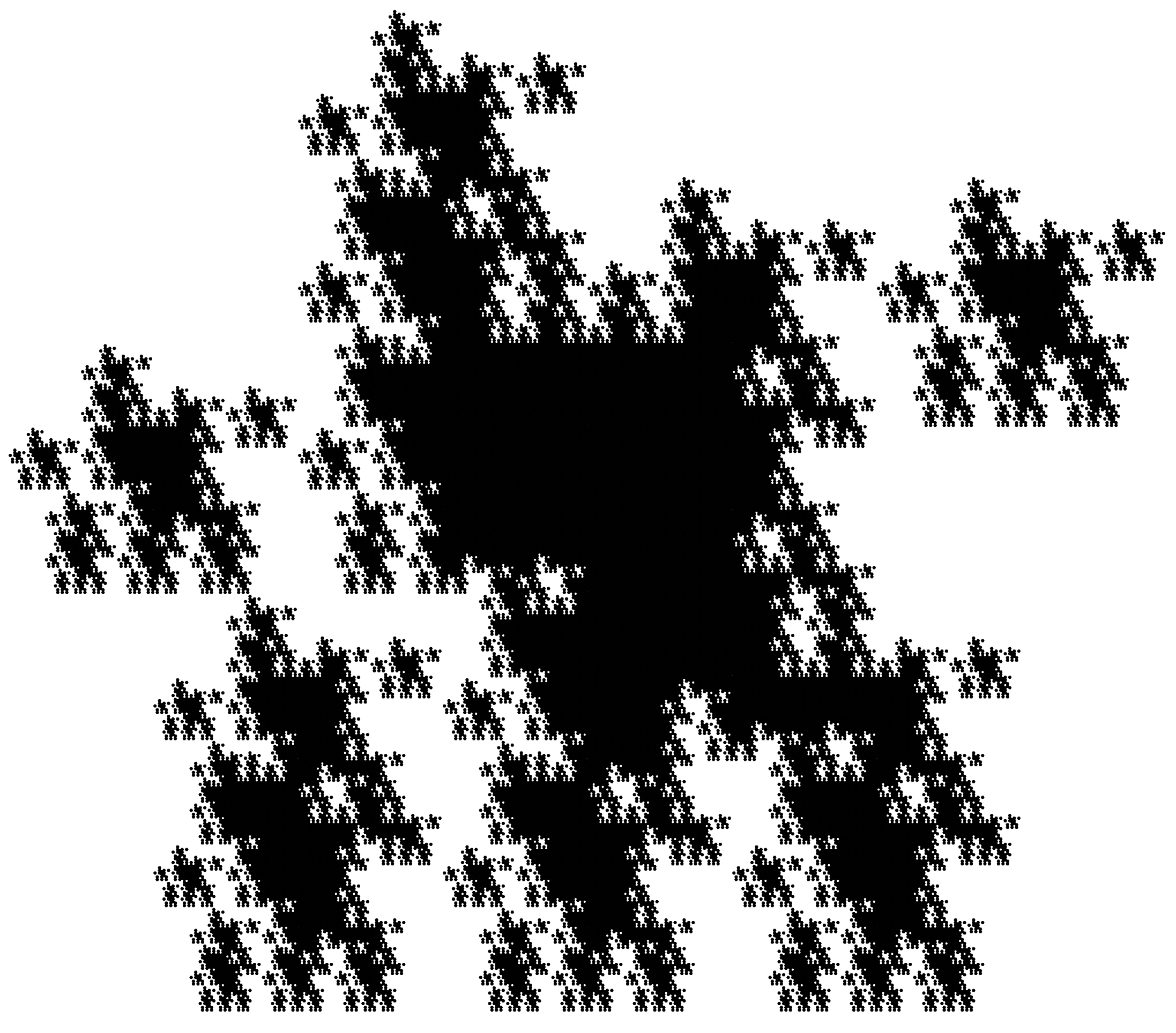}}\qquad}
\caption{\label{fig:strange}
A projected image of a non-layered fractal imaginary cube  ${\mathcal F}(4, D)$.
}
\end{figure}

\begin{proposition}
  A fractal imaginary cube of degree $k$
is projected to a set with positive measure 
along the vector $(ak, bk, 1)$ and its permutations for any integers $a$ and $b$.
\end{proposition}
\begin{proof}
Since the set
$\{(x, y) \mid (x, y, z) \in D\}$ forms the grid $\{0,\ldots,k-1\}^2$ as stated in Lemma \ref{lemmal1}, 
the projected image of $D$ along the vector $(ak, bk, 1)$ to the $xy$-plane, which is
  $\{(x - akz, y - bkz) \mid (x, y, z) \in D\}$, 
  is a complete residue system of $\ZZ^2/k\ZZ^2$.
  Therefore, the statement holds by Corollary \ref{cor-1}.
\end{proof}

Some non-layered fractal imaginary cubes of degree 4 yield projected images 
with positive measures along the vector 
(1,1,1), as illustrated in Figure~\ref{fig:strange}.
However, not all of them exhibit this property, indicating that further research into this area is needed.

\begin{figure}[t]
  \subfloat[The fractal imaginary square $\F^2(3, \D_\OOO)$.]
  {\includegraphics*[height=35mm]{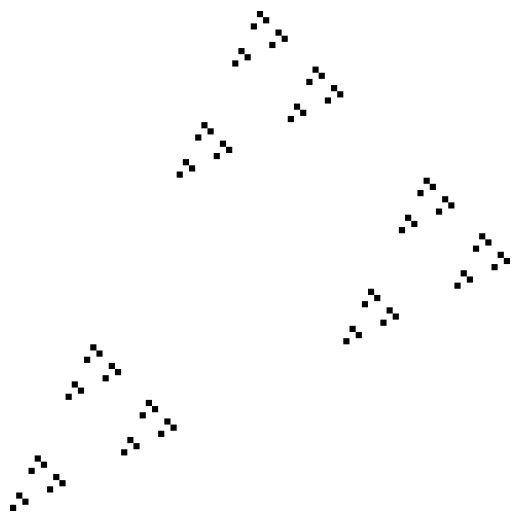}\quad}\qquad
  \subfloat[The one-dimensional Sierpinski gasket $\F^2(3,\D'_\OOO)$.]
  {\quad\includegraphics*[height=35mm]{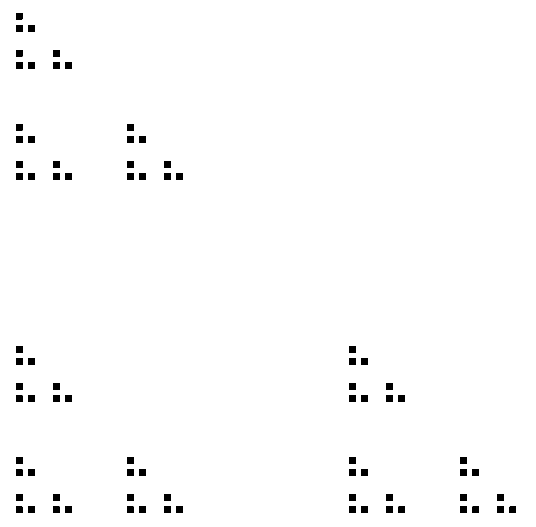}\quad}\qquad
  \subfloat[The differenced radix expansion set $\E^2(3,\Delta(\D'_\OOO))$.]  
  {\includegraphics*[height=35mm]{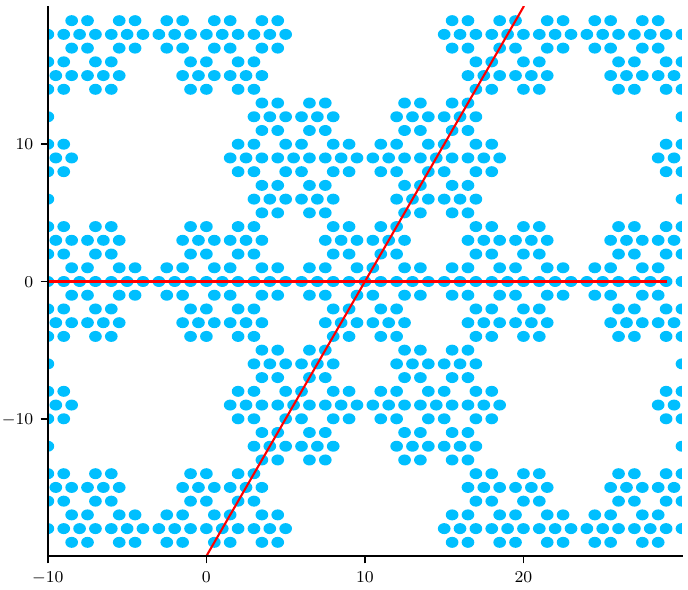}\quad}
   \caption{
The one-dimensional Sierpinski gasket and its associated mathematical structures.
\label{fig:OneDim}} 
 \end{figure}

Another intriguing question would be about fractal imaginary 'squares'.
We define an imaginary square as a two-dimensional object that is projected to  line segments of the same length in orthogonal two directions. A (homothetic) %
 fractal imaginary square of degree $k$ is then defined as an imaginary square that has the form  $\F^2(k, D)$ for $D$ a two-dimensional digit set of cardinality $k$.  A fractal imaginary square of degree $k$ exists corresponding to a permutation of  $\{0,1,\ldots,k-1\}$ following an argument  similar to  Lemma \ref{lemmal1}. 
A line segment is the  only fractal imaginary square of degree 2.
Besides the line segment, the unique fractal imaginary square of degree 3 is $\F^2(3, \D_\OOO)$
for the following digit set $\D_\OOO$  (Figure \ref{fig:OneDim}(a)).
\[
\D_\OOO = \{(0,0), (1,2), (2,1)\}
\]

A layered fractal imaginary square is defined as a fractal imaginary square generated by the set
\[
D^{(2)}_{k, l} = \{(x, y)\in \ZZ^2 \mid 0 \leq x, y < k, \ x + y \equiv l - 1 \bmod k\}
\]
for $k \geq 2$ and $0 \leq l < \frac{k}{2}$.  Accordingly,  $D^{(2)}_{3, 1} = \D_\OOO$, and both $D^{(2)}_{2,0}$ and $D^{(2)}_{3, 0}$ generate line segments.

 The author suggests that
characterizing the directions along which 
layered fractal imaginary squares are projected  
to sets with positive measures may be more challenging
than in the three-dimensional case 
for the following reasons.  Firstly, 
the result similar to Theorem \ref{mainth1} does not hold for the two-dimensional case.
It is straightforward to demonstrate that projected images
of $\F^2(k, D^{(2)}_{k,l})$ along
$(a, b)$ for coprime integers $a, b$ such that 
$a + b$ is coprime to $k$
have positive measures.
This reasoning is analogous to that presented in Lemma \ref{mainlemmapositive}, resulting in  
fractals generated by standard digit sets.
However, the converse is not true.
For instance, consider the two-dimensional digit set $D^{(2)}_{4, 2} = \{(0,1), (1,0), (2, 3), (3,2)\}$. The projection along %
$(3, -5)$ onto the $x$-axis produces the set
$\{0.6, 1, 3.8, 4.2\}$, which becomes $\{0, 1, 8, 9\}$
after an affine transformation. 
This set is an example of a product-form digit set, which is known to 
generate fractals with positive measures (refer to Example 3.1 of \cite  {LagariasWang1}).
Since %
$3 + (- 5) = -2$ 
is a multiple of $2$, 
this example illustrates that direct translation of Theorem 
\ref{mainth1}
 to the two-dimensional case is not feasible.

Secondly, the structure of differenced radix expansion sets for layered fractal imaginary squares significantly diverges
from that observed in the three-dimensional case,
posing challenges to developing a proof akin to Lemma \ref{lemmakl1}.
The fractal produced by
$\D_\OOO$ is an affine image of the one-dimensional 
Sierpinski gasket (\cite{Kenyon1997}), that is $\F^2(3, \D_\OOO')$ for 
\[
\D_\OOO' = \{(0,0), (1, 0), (0, 1)\}
\]
as depicted in Figure \ref{fig:OneDim}(b). 
Kenyon studied the measure-theoretial properties of the projected images of this fractal and established the following:
\begin{theorem}[\cite{Kenyon1997}]\label{th:onedim1}
  $\F^2(3, \D_\OOO')$ is projected to a set with positive Lebesgue measure if and only if the projection is done along a vector $(a, b)$ for coprime integers $a, b$ such that  $a-b$ is a multiple of 3.
\end{theorem}
Through an affine transformation, 
this theorem can be restated for $\F^2(3, \D_\OOO)$ as follows.
\begin{corollary}\label{cor:onedim2}
  $\F^2(3, \D_\OOO)$ is projected to a set with positive Lebesgue measure if and only if the projection is done along a vector $(a, b)$ for coprime integers $a, b$ such that
  $a+b$ is not a multiple of 3.
\end{corollary}
The proof of the equivalence of these two statements is similar to that of the equivalence
of Theorem~\ref{mainth1} and Theorem~\ref{mainlemma}. 
Kenyon proved Theorem \ref{th:onedim1} measure-theoretically through the notion of absolutely continuous invariant measures. 
Let us try to formulate 
a proof of Theorem \ref{th:onedim1} similar to that of Lemma \ref{lemmakl1}.
The differenced digit set of $\D'_\OOO$ is
\[
\Delta(\D'_\OOO) = \{(-1, 1), (-1, 0), (0, -1),(0, 0),(0, 1), (1, 0), (1,-1)\}.
\] 
Note that the fractal generated by $\Delta(\D'_\OOO)$ is the  (affine-transformed) hexaflake fractal, which appeared also as the projectied image of $\HHH_\infty$ along $(1,1,1)$.
The differenced radix expansion set $\E^2(3, \Delta(\D'_\OOO))$ is
depicted in Figure \ref{fig:OneDim}.
For a proof similar to Lemma \ref{lemmakl1}, we need to show the following:
\begin{quote}
  {\textbf Claim.} 
For every pair of coprime
integers
$(a,b) \in \ZZ^2$ 
such that
$a - b$ is not a multiple of 3, there exists %
$j \ne 0$ 
such that
$(ja, jb) \in \E^2(3, \Delta(\D'_\OOO))$.
\end{quote}
Nevertheless, direct proof of this claim appears to be challenging.
Specifically, a strategy analogous to that used in Lemma \ref{lemmakl1} may not be feasible, because of the numerous 'holes' that are present in $\E^2(3, \Delta(\D'_\OOO))_{3^m(3^n-1)}$.
We would like to list the above claim as a corollary to Theorem \ref{th:onedim1}, because 
Theorem \ref{theorem:delta} can be generalized to projections of $n$-dimensional digit sets, and Theorem \ref{theorem:delta} claims an equivalence. 

Kenyon also explored measures and Hausdorff dimensions of projections of the
one-dimensional Sierpinski gasket. Moreover,
Hochman demonstrated in \cite{Hochman} that the Hausdorff dimensions of irrational
projections of the one-dimensional Sierpinski gasket are zero.
As outlined in \cite{ProjectionBook1} and
\cite{ProjectionBook2}, 
the study of Hausdorff dimensions of projected images of fractals is considered an intriguing area of research. 
The investigation of Hausdorff dimensions of the projected images of fractal imaginary cubes remains an open field of inquiry.
\bigskip

{\bf Acknowledgement. }
  The author would like to thank Takayuki Tokuhisa for illuminating discussions.
  This work was supported by JSPS KAKENHI Grant Number 15K00015, 23K28036, and 18H03203.

\bibliographystyle{emss}
\bibliography{ref1}

\newpage
\appendix

\section{Projections of 3D-printed models of fractal imaginary cubes}

Captions specify projection directions under the identification of the objects with $\SSS_\infty$, $\TTT_\infty$, and $\HHH_\infty$. Red frames indicate that the shadows have positive measure even when the
models were perfect mathematical representations.

\medbreak

\noindent
\begin{tabular}{ccccc}
    \multicolumn{5}{c}{Sierpinski Tetrahedron}\\
    \includegraphics[height=3.5cm,width=2.6cm]{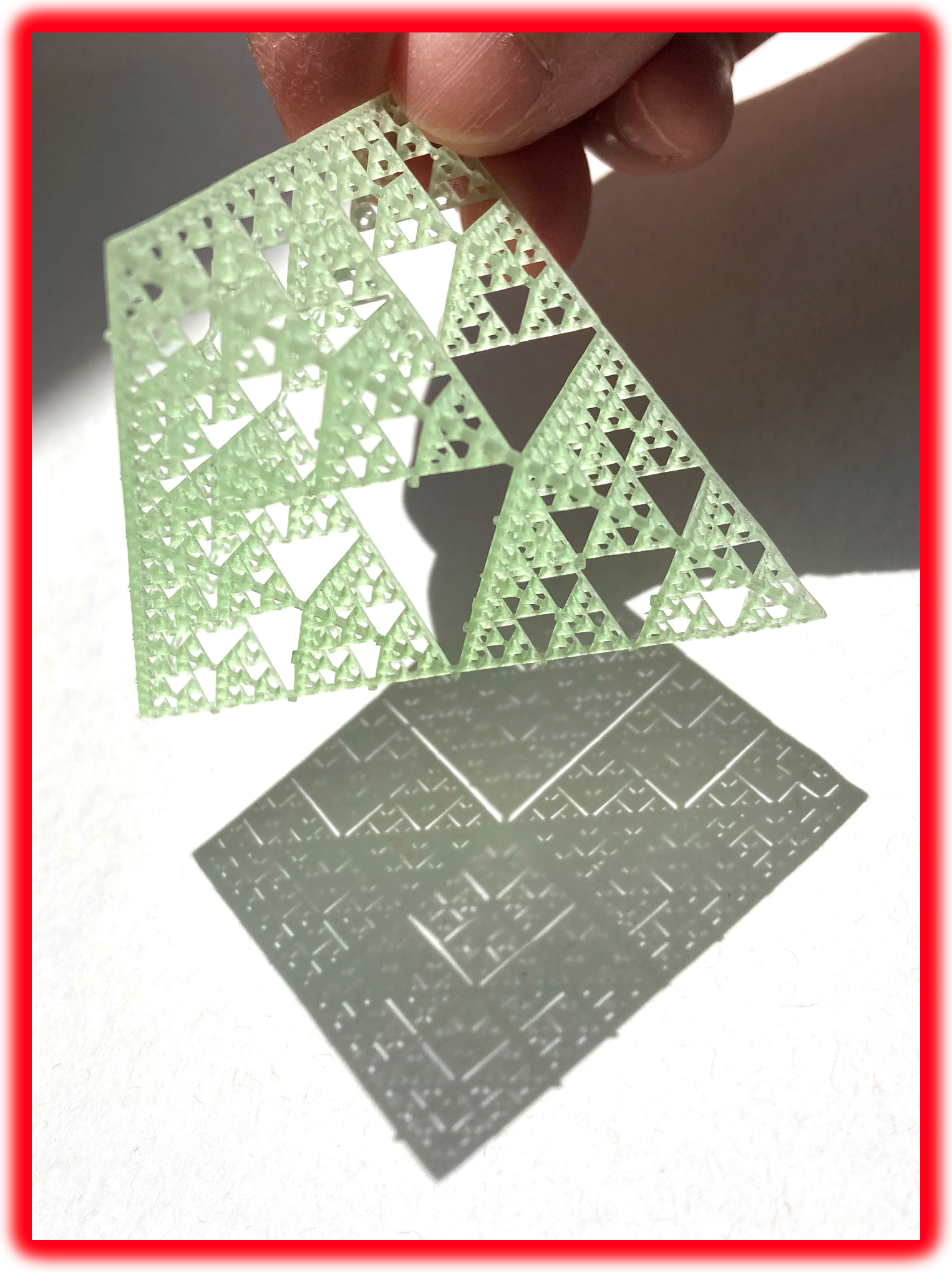} &
\includegraphics[height=3.5cm,width=2.6cm]{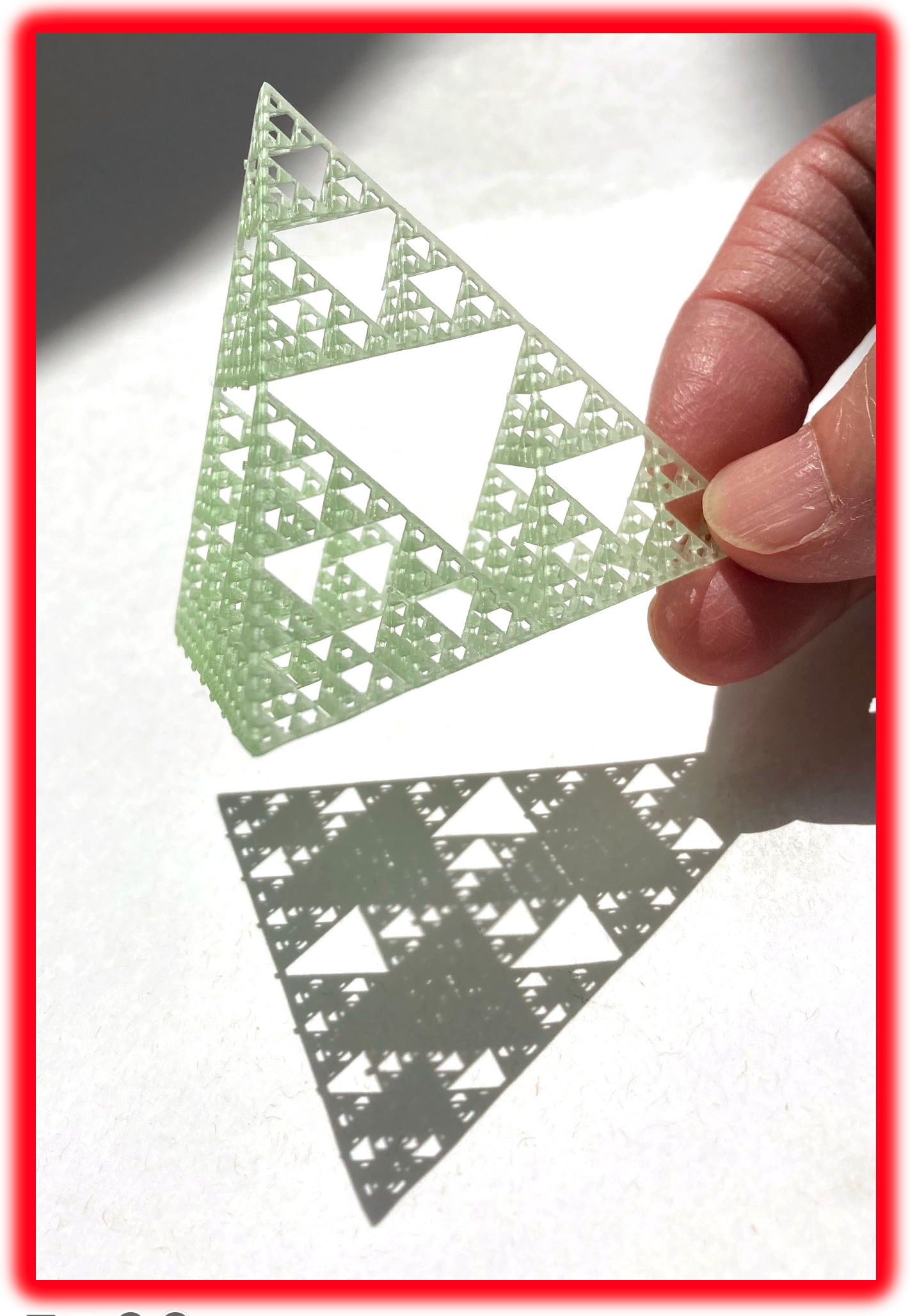}  &
\includegraphics[height=3.5cm,width=2.6cm]{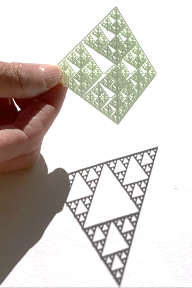}  &
\includegraphics[height=3.5cm,width=2.6cm]{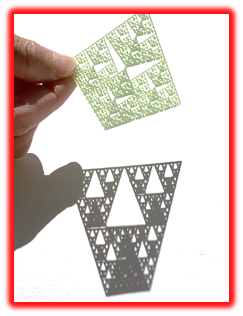}  &
\includegraphics[height=3.5cm,width=2.6cm]{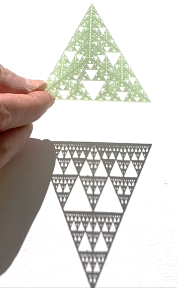} \\
(0,0,1)& (1,1,1)& (1,1,0) & (1,0,2) &(1,1,-2)\\\\
\multicolumn{5}{c}{T-Fractal}\\
\includegraphics[height=3.5cm,width=2.6cm]{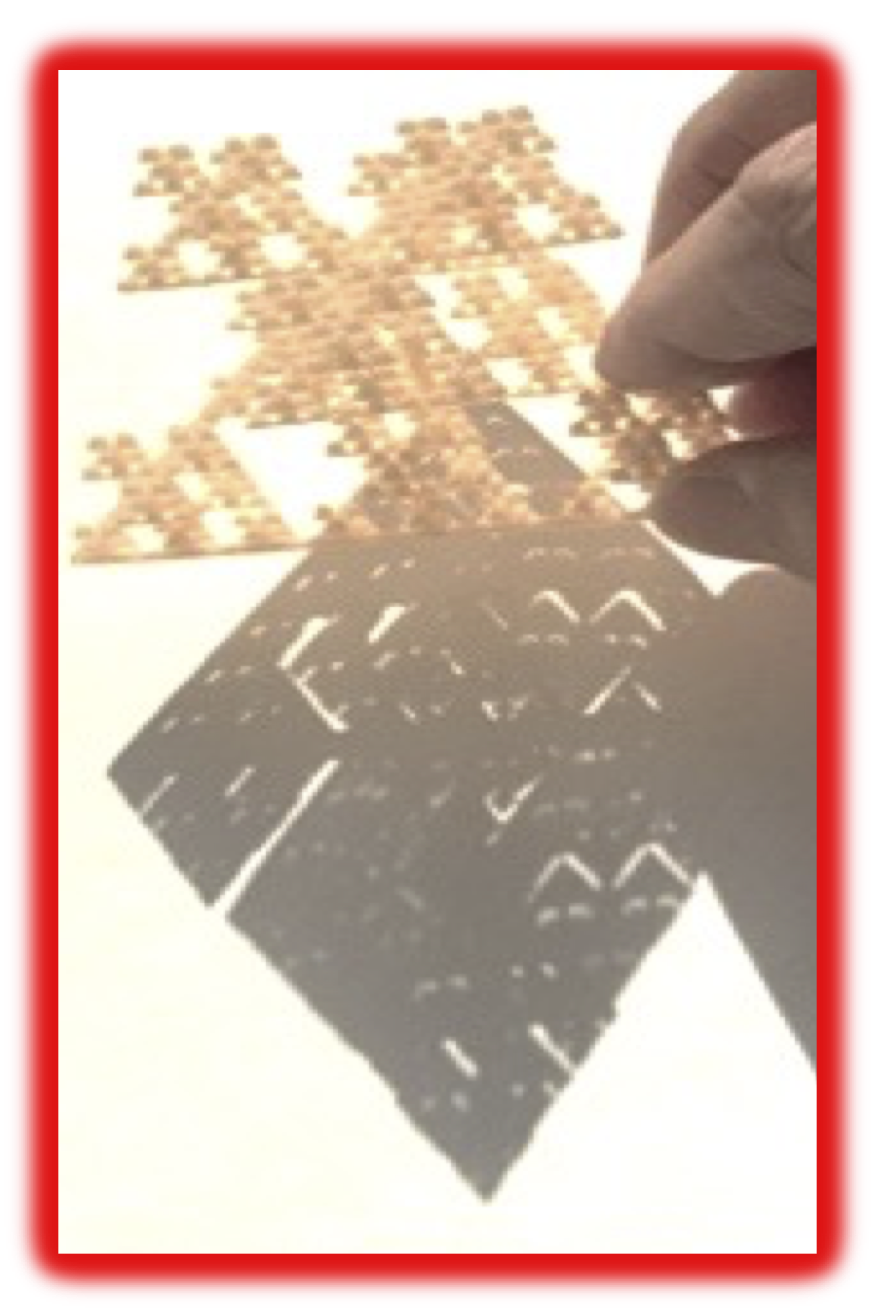}  &
\includegraphics[height=3.5cm,width=2.6cm]{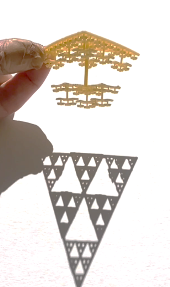}  &
\includegraphics[height=3.5cm,width=2.6cm]{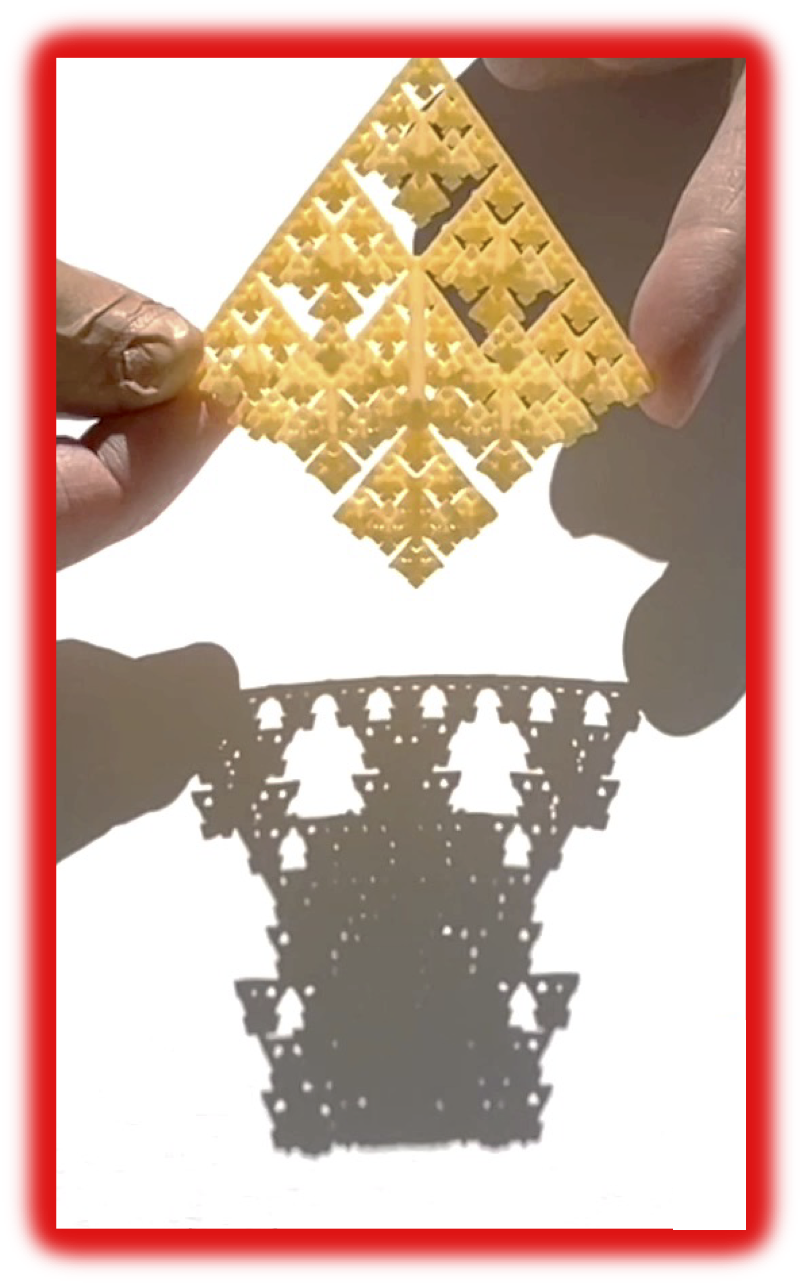}  &
\includegraphics[height=3.5cm,width=2.6cm]{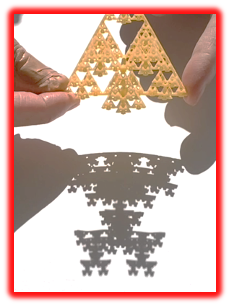}  &
\includegraphics[height=3.5cm,width=2.6cm]{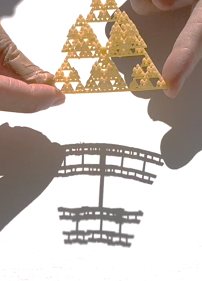} \\
(0,0,1)& (1,1,1)& (1,1,0) & (1,1,-1) &(1,1,-2)\\\\

\multicolumn{5}{c}{H-Fractal}\\
\includegraphics[height=3.5cm,width=2.6cm]{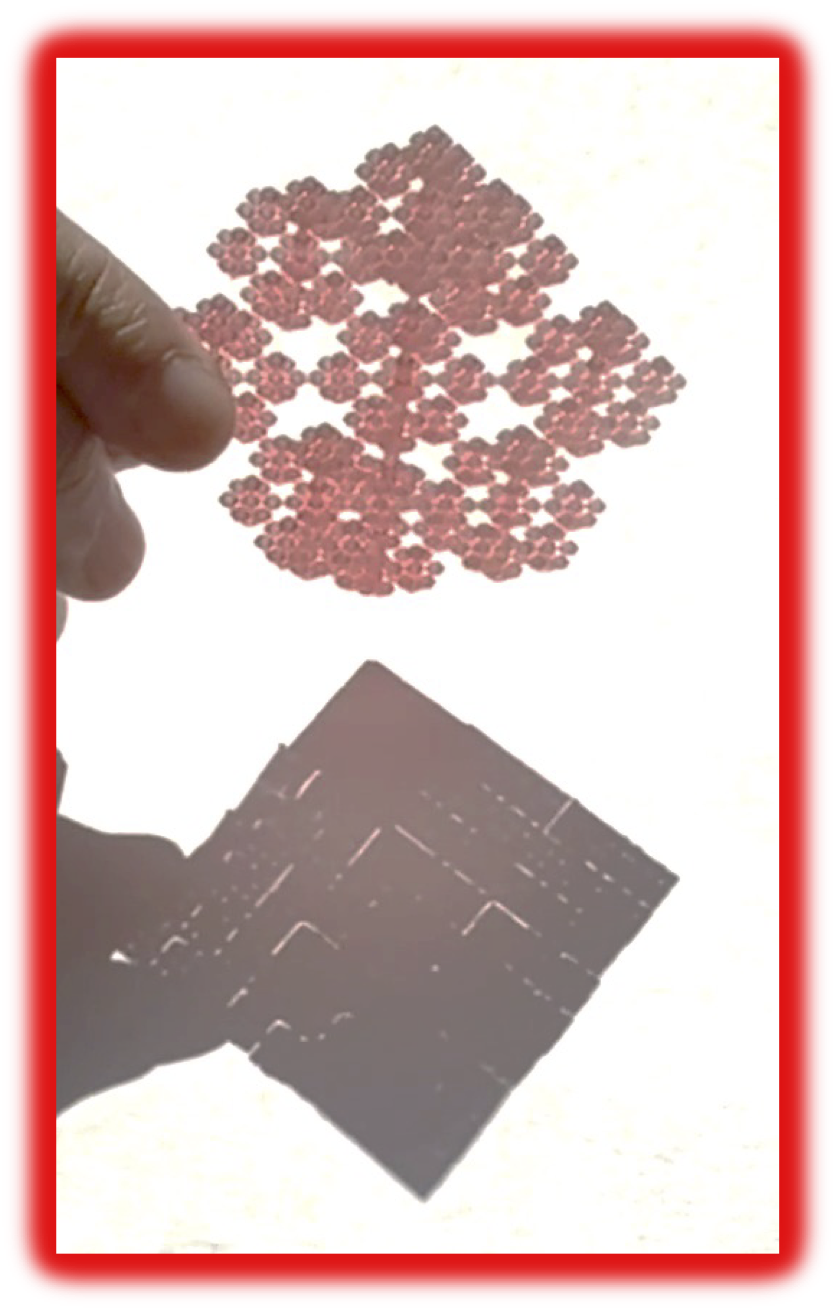}  &
\includegraphics[height=3.5cm,width=2.6cm]{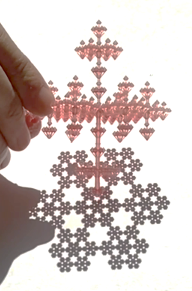}  &
\includegraphics[height=3.5cm,width=2.6cm]{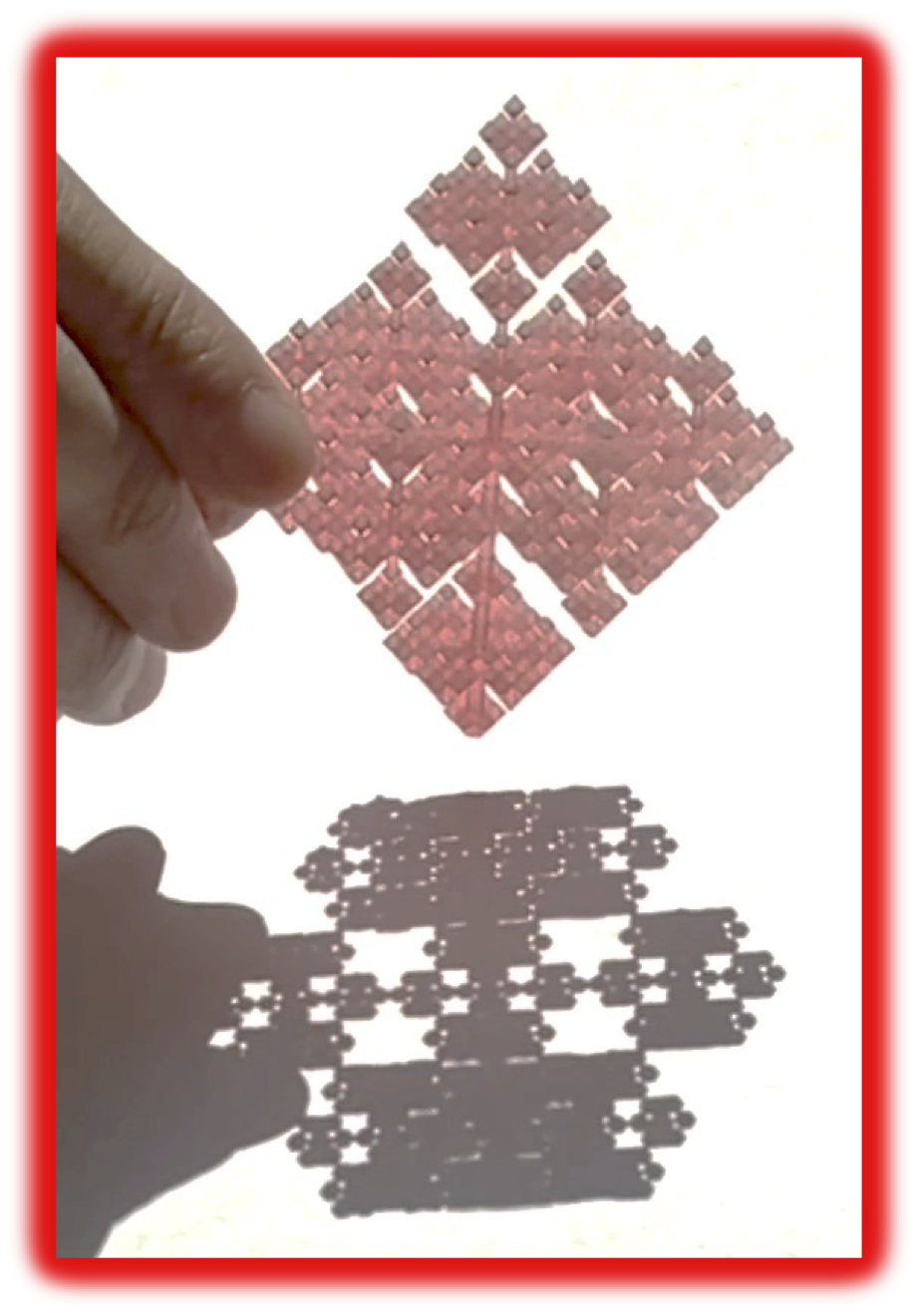}  &
\includegraphics[height=3.5cm,width=2.6cm]{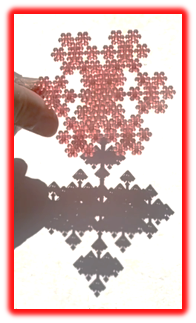}  &
\includegraphics[height=3.5cm,width=2.6cm]{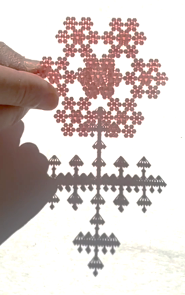}  \\
(0,0,1)& (1,1,1)& (1,1,0) & (1,1,-1) &(1,1,-2)\\
(2,2,1) &        & (1,1,4) & (1,1,-5) \\
\end{tabular}

\end{document}